\pdfoutput=1
\documentclass[ a4paper 
              , UKenglish
              , cleveref
              , autoref
              , thm-restate
              ]{lipics-v2021}

\hideLIPIcs

\usepackage{fouche}
\usepackage{tangle}

\def\oneCat#1{\mathsf{#1}}
\let\cate\oneCat
\def\twoCat#1{\underline{\oneCat{#1}}}
\def\dblCat#1{\mathbf{#1}}
\def\dblD{\dblCat{D}}
\def\dblC{\dblCat{C}}
\def\proc#1{\twoCat{KSW}(#1)}
\def\pro{\rightsquigarrow}
\let\to\longrightarrow
\def\Loose{\mathrel{\ooalign{$\Rightarrow$\cr$\hfil\raise0.1ex\hbox{\rule{0.4pt}{1ex}}\kern3pt\hfil$}}}
\NewDocumentCommand{\smallCell}{ m m m m O{} }{%
	\raisebox{-1pt}{\begin{tikzpicture}[baseline=(cell.south)]
		\node[font=\small, draw, minimum height=0em, minimum width=1em, inner sep=2pt] (cell) {$#5$};
		\node[font=\tiny, above=0em of cell, inner sep=1pt] {$#1$};
		\node[font=\tiny, left =0em of cell, inner sep=1pt] {$#2$};
		\node[font=\tiny, right=0em of cell, inner sep=1pt] {$#3$};
		\node[font=\tiny, below=0em of cell, inner sep=1pt] {$#4$};
\end{tikzpicture}}%
}
\NewDocumentCommand{\dblMly}{ O{\clK} }{\hyperref[def_dbl_processes]{\dblCat{KSW} (#1)}}
\def\cartDblMly{\hyperref[cart_dbl_mly]{\dblCat{Mly}(\clK)}}

\def\emptyList{[\,]}
\def\cons{\mathbin{::}}

\makeatletter
\newcommand{\concat}{\mathbin{\mathpalette\concat@inner\relax}}
\newcommand{\concat@inner}[2]{%
	\ooalign{$\m@th#1\scriptscriptstyle\mathplus$\cr\hfil$\m@th#1\mkern6mu\mathplus$\hfil\cr}%
}
\newcommand{\mathplus}{\mathord{\mkern2mu\raisebox{0.15ex}{\scalebox{0.85}{\texttt{+}}}\mkern2mu}}
\makeatother

\def\idpto#1{\begin{tikzcd}[ampersand replacement=\&, row sep=tiny,
			nodes={inner sep=1pt, outer sep=2pt}] {#1} \ar[r, "\shortmid"{marking}, Rightarrow, no head] \& {#1} \end{tikzcd}}

\usepackage{string-diagrams}

\def\dact{\otimes^+}
\def\sact{\odot^+}
\def\cotabu{\intercal^*}
\def\tabu{\intercal_*}
\def\inj{\text{inj}}
\def\vertC{\odot_\text{v}}
\def\horC{\odot_\text{h}}

\definecolor{dColor}{HTML}{006400}
\definecolor{sColor}{HTML}{F5DEB3}
\definecolor{muColor}{HTML}{A52A2A}

\newtheorem{notation}[theorem]{Notation}

\usepackage{xspace}
\def\PDC{pseudo double category\@\xspace}
\usepackage{quiver}
\usepackage{environ} 
\newif\ifshowproofs
\showproofstrue 

\NewEnviron{savvyProof}[1][Proof]{%
  \ifshowproofs
    \begin{proof}[#1]
    \BODY
    \end{proof}
  \fi
}

\title{Monads and limits in bicategories of circuits}
\author%
{Fosco Loregian}
{Tallinn University of Technology, Tallinn, Estonia}
{fosco.loregian@taltech.ee}
{}
{}

\authorrunning{F. Loregian}

\Copyright{Jane Open Access and Joan R. Public} 

\ccsdesc[500]{Theory of computation~Formal languages and automata theory}

\keywords{Bicategory, double category, monad, cotabulator, proarrow equipment.}

\begin{document}
\maketitle
\begin{abstract}
	We study monads in the (pseudo-)double category $\dblMly$ where loose arrows are \emph{Mealy automata} valued in an ambient monoidal category $\clK$, and the category of tight arrows is $\clK$. Such monads turn out to be elegantly described through instances of \emph{semifree bicrossed products} (bicrossed products \emph{of monoids}, in the sense of Zappa-Szép-Takeuchi, where one factor is a free monoid). This result which gives an explicit description of the `free monad' double left adjoint to the forgetful functor. (Loose) monad maps are interesting as well, and relate to already known structures in automata theory. In parallel, we outline what double co/limits exist in $\dblMly$ and express in a synthetic language, based on double category theory, the bicategorical features of Katis-Sabadini-Walters `bicategory of circuits'.
\end{abstract}
\section{Introduction}
\begin{definition}\label{def_zero}
	Let $\twoCat{B}$ be a bicategory, and regard the monoid $\bbN$ of natural numbers as a bicategory with a single object and only identity 2-cells $\id_n : n\To n$; in \cite{Katis1997,Katis2002} the authors construct the bicategory $\Omega\twoCat{B}$:
	\begin{itemize}
		\item 0-cells are the \emph{pseudofunctors} $F : \bbN\to \twoCat{B}$; so, a pair $(X,a : X\to X)$ of an object and an endo-1-cell of $\twoCat{B}$.
		\item 1-cells $\alpha : F \To_l G$ are the \emph{lax natural transformations}, consisting of intertwiners $(f,\varphi) : (X,a) \to (Y,b)$ where $f : X\to Y$ is a 1-cell and $\varphi : f \circ a \To b\circ f$ is a 2-cell.
		\item 2-cells $\Theta : \alpha \Tto \beta$ are the \emph{modifications} of such lax natural transformations.
	\end{itemize}
	When $\twoCat{B}=\Sigma\clK$ is the one-object bicategory naturally associated to a monoidal category $(\clK,\otimes,I)$, the bicategory $\proc\clK=\Omega\Sigma\clK$ of \emph{processes} has
	\begin{enumtag}{p}
		\item 0-cells the object $A : \star\to \star$ of $\clK$;
		\item 1-cells $A\to B$ the pairs $(E,\varphi)$ where $E$ is an object in $\clK$, and $\varphi : A\otimes E\to E\otimes B$;
		\item \label{ducelles} 2-cells $(E,\varphi) \To (E',\psi)$ the morphisms $\alpha : E\to E'$ in $\clK$ such that $\psi\circ (A\otimes\alpha) = (\alpha\otimes B)\circ \varphi$.
	\end{enumtag}
\end{definition}
\begin{remark}\label{florpo}
	The interest for such an object is easily seen: in the case where $\clK$ is a Cartesian category, as a consequence of the universal property of Cartesian products, a 1-cell in $\proc\clK$ splits as a span
	$\footnotesize\begin{tikzcd}[cramped]
			E & {A\times E} & B.
			\arrow[from=1-2, to=1-1, "d"']
			\arrow[from=1-2, to=1-3, "s"]
		\end{tikzcd}$
	i.e. as a \emph{Mealy automaton} \cite{Ehrig} where $A,B$ are input and output objects, $d$ is an action of $A$ on $E$ and $s$ a morphism yielding a \emph{final state}. In the general monoidal case, a 1-cell is a process $A\otimes E \to E\otimes B$ of a similar shape, where the output and the final state can correlate.
\end{remark}
Rather surprisingly, the properties of $\proc\clK$ \emph{qua} bicategory (i.e. its exactness or co/completeness properties, its formal theory of monads, etc.) have not been the subject of a thorough analysis, not even in openly category-theoretic oriented derivative works of \cite{Katis1997} focused on representing categories of processes \cite{DeFrancescoAlbasini2011,Katis2002,ROSEBRUGH_SABADINI_WALTERS_1998}, finding their universal properties \cite{10.1145/3531130.3533365,DiLavore2023}, exploiting their nature of two-dimensional objects or generalising their features \cite{bonchi2024effectfulmealymachinesbisimulation}; albeit enough for most applications, truncations or simplifications of the bicategory in study are fairly unnatural from a category-theoretic standpoint.

The present work is meant to bridge this gap and fits into a line of works \cite{boccali2023semibicategory,EPTCS397.1,boccali_et_al:LIPIcs.CALCO.2023.20,loregian2024automata} with a similar goal of using 2-dimensional category theory as a foundation for the theory of automata (while at the same time, expanding and contributing to the state of the art of the former discipline).

As a starting point, we embed the bicategory $\proc\clK$ into the far more natural\footnote{An example of why this should be a more natural structure to consider, upon which we will expand later on: \cite{Katis1997} claims that $\proc\clK$ has the structure of a \emph{symmetric monoidal bicategory}, clearly underestimating how daunting the verification of \emph{all} axioms is --following \cite{coherence-tricat}, a monoidal bicategory is a tricategory with a single 0-cell, a gadget the definition of which runs to many pages. At the time of writing, the only slick way to generate a monoidal structure on a bicategory $\twoCat{B}$ is to present $\twoCat{B}$ as the bicategory of loose arrows of a \PDC, and verify a certain technical condition, cf. \cite[Definition 4.10]{hansen2019constructingsymmetricmonoidalbicategories}.} structure of a \emph{(pseudo) double category} (\cite{CTGDC_1999__40_3_162_0,Grandis2019,Shulman2008a}), whose cells have the form
$\footnotesize\begin{tikzcd}[cramped]
		A & B \\
		X & Y
		\arrow["{\langle E,d,s\rangle}", "\shortmid"{marking}, from=1-1, to=1-2]
		\arrow["f"', from=1-1, to=2-1]
		\arrow["\alpha"{description}, draw=none, from=1-1, to=2-2]
		\arrow["g", from=1-2, to=2-2]
		\arrow["{\langle E,d',s'\rangle}"', "\shortmid"{marking}, from=2-1, to=2-2]
	\end{tikzcd}$
thus having $\clK$ as its 1-category of tight arrows, and $\proc\clK$ as bicategory of loose arrows.

The theory of monads, comonads, limits, colimits, Kan extensions, \dots can now be re-enacted in this \PDC: we find how some co/limits of a double categorical nature exist, and some do not: a specific obstruction that prevents many colimits from existing appears in \autoref{gdaogao}.

Most of our interest lies in the Cartesian case,\footnote{In fact, for the purposes that we set in the present paper of characterising monads in the \PDC of Mealy automata (cf. \autoref{sgrulbio}) and their left modules (cf. \autoref{crucial_thm}), we additionally simplify things pretending that computations happen in $\Set$, but this is a minor issue as `element-free' proofs of the various statements can easily be found.} where the representation in \autoref{florpo} for 1-cells of $\proc\clK$ is available; in this case we call $\cartDblMly$ the double category in study, to reflect that it reduces to a double category `of Mealy automata'. A similar study of monads and limits can be conducted for a generic monoidal category $\clK$ (an important example of which is the Kleisli category of a strong or commutative monad on $\Set$, like the finite distributions monad \cite{Giry1982,Fritz2020,Petri_an_2021,Jacobs2014}, or the Maybe monad $\firstblank+1$, or\dots), and some of its ideas lie at the heart of a fairly abstract branch of representation theory \cite{Majid1991,Kassel_1995}.

The restriction to the Cartesian case has multiple motivations: one has first to resolve the simplest case of the more general taxonomy problem of a double category defined out of \cite{Katis1997}'s $\Omega\twoCat{B}$ (which seems to be buried in the aforementioned literature on Hopf algebras). Second, the  Cartesian case is already interesting from the category-theoretic standpoint (but we do not feign that the endeavour here has any motivation beyond finding application and interesting examples for the language of double category theory).

In \autoref{bonbiollo} we superficially sketch how things might change in a generic monoidal category, with special attention to the core results of the paper, \autoref{crucial_constru}, \autoref{crucial_thm}.

A short summary of the paper is easily given.
\begin{itemize}
	\item (\autoref{sgrulbio}, \autoref{dgaddaaa}) A monad $M : A\pto A$ in $\cartDblMly$ is determined by a pair $(d_M,s_M)$ of compatible monoid actions, one of which is the free monoid on $A$, and the other is the state space $E$ of the automaton $E \xot{d_M} A\times E \xto{s_M} A$. As a consequence, \autoref{coro_fremonad} allows to describe free monads on loose endomorphisms $F : A\pto A$.
	\item (\autoref{crucial_thm}) The category of algebras for a monad $M$ as above coincides with the category of left representations for the \emph{bicrossed product} \cite{Agore2009forfinite} $E \bowtie A^*$ of monoids obtained from the compatible actions above.
	\item (\autoref{horatio}) Every adjunction gives rise to a monad, but monads are way more than those obtained from adjoint pairs, as the constraints for a pair of loose morphisms in $\cartDblMly$ to be adjoint are pretty strict. This rules out the existence of an Eilenberg-Moore construction for monads in $\cartDblMly$.
	\item If $\clK$ is co/complete, $\cartDblMly$ has all `local co/limits', but it lacks most interesting global co/limits; it has a global terminal object; it lacks (for trivial reasons) an initial one (a fact that we observe very early, \autoref{gdaogao}), as well as coproducts of any cardinality, \autoref{coprod_no}; however, it has at least one kind of genuinely double-categorical colimit, cotabulators (cf. \autoref{constr_cotab}).
\end{itemize}
The remainder of the section builds (in the general monoidal case) the \PDC that is the subject of this paper. We recall the fundamentals needed for the theory of double categories, but since displaying all axioms is space-consuming and hinders readability, we relegate almost all detailed definitions of double category theory to the appendix. \cite[§7, I.3.3]{CTGDC_1999__40_3_162_0} is a classic reference from which we draw the basic definitions. One can also consult \cite{niefield2012span} (also a precious source of examples, see [\emph{ibi}, 2.1---2.8]).

The theory of (strict) double categories is very easy to summarise in a nutshell: the essentially algebraic theory \cite[3.D]{Adamek1994} of categories is interpreted in the category $\Cat$ of categories and functors; a double category is thus just a model of such a theory, i.e. a category $\dblD=(\clD_0,\clD_1)$ internal to $\Cat$, thus equipped with composition, source, target and identity functors.\footnote{The reader in need of a more precise definition can refer to \cite[Chapter 8]{Bor1}, where internal category theory in a category $\clE$ with finite limits is reviewed.}

Many examples of such a structure however arise when the category of objects $\clD_0$ does not have a strictly associative and unital composition law, but where these properties only hold up to coherent isomorphisms living in $\clD_1$; this makes the composition in one direction (loose) weak, while the composition in the other direction (tight) remains strict.

The focus on strict double categories however isn't a true restriction; similarly to what happens for bicategories, every \PDC $\dblD$ can be strictified into a double category $\dblD^\parallel$ where the loose morphisms compose on the nose. Such strictification procedures, although of the highest conceptual importance, are not suitable to take as a working definition.
\subsection{The \PDC \texorpdfstring{$\dblMly$}{Mly}}
The definition of \PDC is recalled below in \autoref{def_pdc}.

Fix a monoidal category $(\clK,\otimes,I)$; the subject of the present work will be the \PDC obtained combining $\clK$ (as tight or `vertical' category $\clV(\dblMly)$) and Walters' $\proc\clK$ (as loose, or `horizontal' bicategory $\clH(\dblMly)$), as follow.
\begin{definition}\label{def_dbl_processes}
	The \emph{\PDC of processes} $\dblMly$ has
	\begin{itemize}
		\item objects the same of $\clK$, denoted as $A,B,C,\dots,X,Y,Z,\dots$ and the like;
		\item $\clK$ as category of tight morphisms, denoted as $f,g,h,\dots,u,v,w,\dots$ and the like;
		\item the bicategory $\proc\clK$ of \autoref{def_zero} as bicategory of loose morphisms; we denote loose morphisms as $\langle E,\phi\rangle : A \pto B$ where $\phi : A\otimes E\to E\otimes B$ is a morphism in $\clK$;
		\item choice of cells determined as follows:
		      the set of cells having
		      $\footnotesize\begin{tikzcd}[cramped]
				      A & B \\
				      X & Y
				      \arrow[""{name=0, anchor=center, inner sep=0}, "{\langle E,\phi\rangle}", "\shortmid"{marking}, from=1-1, to=1-2]
				      \arrow["f"', from=1-1, to=2-1]
				      \arrow["g", from=1-2, to=2-2]
				      \arrow[""{name=1, anchor=center, inner sep=0}, "{\langle E',\psi\rangle}"', "\shortmid"{marking}, from=2-1, to=2-2]
			      \end{tikzcd}$
		      as boundary consists of the morphisms $\alpha : E\to E'$ in $\clK$ such that the square
		      \[\label{baiadiai}\begin{tikzcd}[cramped]
				      {A\otimes E} && {E\otimes B} \\
				      {X\otimes E'} && {E'\otimes Y}
				      \arrow["\phi", from=1-1, to=1-3]
				      \arrow["{f\otimes\alpha}"', from=1-1, to=2-1]
				      \arrow["{\alpha\otimes g}", from=1-3, to=2-3]
				      \arrow["\psi"', from=2-1, to=2-3]
			      \end{tikzcd}\]
		      commutes (evidently, this is but  slight modification of \autoref{ducelles}, which recovers (and thus extends nontrivially) $\proc\clK$ as the double subcategory $\clH(\dblMly)$ of \emph{tightly discrete} cells --where $f,g$ are both identities of the respective objects).
	\end{itemize}
\end{definition}
\begin{remark}\label{gdaogao}
	A fundamental observation for what follows is that the bicategory of those tightly discrete cells and the whole $\dblMly$ behave in a radically different fashion with respect to their co/completeness properties: the existence of a cell $\alpha$ as in \eqref{baiadiai} is governed by the existence of morphism between carriers of the processes $\phi,\psi$, and thus there can be few of them. For example, no matter how cocomplete $\clK$ is, $\dblMly$ lacks an initial object, which (cf. \cite{CTGDC_1999__40_3_162_0}) would consist of a unique cell
	$\footnotesize\begin{tikzcd}[cramped]
			U \otimes I \ar[d, "!\otimes\upsilon"']\ar[r,"\shortmid"{marking}] & I \otimes U \ar[d, "\upsilon\otimes!"]\\
			A \otimes E \ar[r,"\shortmid"{marking}, "\varphi"'] & E \otimes B
		\end{tikzcd}$
	for every $(E,\phi) : A \pto B$; the only choice to have a unique $U\to A,U\to B$ for every $A,B$ is to take $U=\varnothing$ (an initial object in $\clK$), but then, we remain with the problem that $\upsilon : I\to E$ exists (the same for every $E$) even though the loose domain on $\upsilon$ `reduces' to the isomorphisms $\varnothing\otimes I \cong I\otimes \varnothing\cong \varnothing$.

	On the other hand, $\proc\clK$ admits \emph{local} initial objects, in the sense that each category $\proc\clK(A,B)$ has an initial object as soon as $\clK$ does (and way more: it is cocomplete, or locally presentable, as soon as $\clK$ is cocomplete or locally presentable, as a consequence of the pullback characterisation in \cite[Proposition 3.5]{boccali_et_al:LIPIcs.CALCO.2023.20}).
\end{remark}
In fact, this is an instance of a much more general obstruction (for example, it prevents binary coproducts from existing).
\begin{remark}\label{a_trick}\label{esem_cell}
	Some universal objects in $\dblMly$ have the form
	$
		\begin{tikzcd}[cramped]
			X & X \\
			A & B
			\arrow["\shortmid"{marking}, Rightarrow, no head, from=1-1, to=1-2]
			\arrow["f"', from=1-1, to=2-1]
			\arrow["g", from=1-2, to=2-2]
			\arrow["{\langle E,\phi\rangle}"', "\shortmid"{marking}, from=2-1, to=2-2]
		\end{tikzcd}$
	i.e. they have as loose domain the identity cell $\iota_X : X\pto X$, which can be easily seen to consist of the isomorphism $\langle I,\lambda_X^{-1}\circ\rho_X\rangle : X\otimes I \to I\otimes X$ induced by the unitors of $(\clK,\otimes)$.

	In particular, when $\clK$ is Cartesian, this implies that a cell $\alpha : \top\to E$ like \eqref{esem_cell} points $E$ in a certain canonical way (in particular, the state space $E$ can't be empty), and that (when $\phi =\langle d,s\rangle$ splits into a span as in \autoref{florpo}) the following constraints for the existence of a cell like \eqref{esem_cell} must be satisfied:
	\begin{itemize}
		\item compatibility with $d$ (notation as in \autoref{florpo}): for all $x\in X$, one has $d (f x , \alpha(\top)) = \alpha(\top)$;
		\item compatibility with $s$ (notation as in \autoref{florpo}): for all $x\in X$, one has $s (f x , \alpha(\top)) = g x$.
	\end{itemize}
	Depending on the nature of $d,s$ these requirements can be met or not; this remark will turn out particularly useful in \autoref{genuinely} to rule out the existence of some co/limits in $\dblMly$.
\end{remark}
\begin{notation}\label{cart_dbl_mly}
	The main results of the paper will assume that $\clK$ is Cartesian, as general monoidal structures would not allow a characterisation of monads and algebra as clean as our \autoref{sgrulbio} and \autoref{crucial_thm}; however, we will suggest how some results might be exported (cf. \autoref{bonbiollo}).

	In case $\clK$ is Cartesian monoidal, and thus we can represent loose morphisms of $\proc{\clK,\times}$ as spans like \autoref{florpo}, such loose morphisms $\fkf : A\pto B$ are Mealy automata $(E,\langle d,s\rangle) : A\times E \to E\times B$; we will refer to $\dblMly[\clK,\times]$ as the \emph{double category of} ($\clK$-valued) \emph{Mealy automata} and denote it $\cartDblMly$.
\end{notation}

\section{The structure of \texorpdfstring{$\cartDblMly$}{Mly}}
The core result of the section (and of the paper) characterizes monads in the double category $\cartDblMly$ (and as a consequence, in its loose bicategory $\proc{\clK,\times}$); the result (\autoref{sgrulbio}, \autoref{crucial_thm}) is stated in the category of sets and functions just for convenience: proofs can easily be made diagrammatical, so as to adapt to the case of a general Cartesian category (admitting the co/limits that are needed, for example, to build $A^*$ as a countable coproduct); the full generalization to a monoidal category $(\clK,\otimes,I)$ is also possible, but for a different model of automata in a monoidal category, e.g. that of \cite{Ehrig}.\footnote{In that case, one has to replace the parts of the statement with the appropriate notion of (Hopf) monoid in \autoref{sgrulbio} or module \autoref{crucial_thm} --bicrossed products exist, for example, in the category of Hopf algebras over a field $k$. See \autoref{how_translate} for a more detailed analysis. However, there is no double category where monoidal automata understood as spans $E\xot d A\otimes E \xto s B$ are 1-cells, if $\otimes$ is not the Cartesian structure.}

The structure of the section is quite simple: each subsection introduces a certain double categorical concept; then we analyze its existence (or non-existence) and shape in the particular case of $\cartDblMly$; sometimes there is something to learn from a non-existence proof: for example, in some cases the obstruction is intrinsic, and in other it is a consequence of the excessive strictness of the double categorical concept in question. Then, we draw consequences from the result in study linking it with the existing literature on similar matters if any.

\subsection{Monoidality}
The aim of this subsection is to close a gap in how the \PDC of Mealy automata is (rightly) regarded as a symmetric monoidal bicategory, with no proof but a claim in \cite{Katis1997} (with the definition of SMB appearing in print just a year before, in \cite{carmody1995cobordism}); constructing SMBs is a well-known daunting task, with a panoply of conditions that are often burdensome to check, as they consist of asserting the commutativity of notoriously complicated 3-dimensional diagrams, \cite{coherence-tricat,stay2013compact}.

\cite[1.1]{hansen2019constructingsymmetricmonoidalbicategories} provides instead a criterion to check that the loose bicategory of a \PDC is an SMB.
\begin{proposition}\label{shllo_rula}
	If $\dblD$ is a monoidal double category, of which the monoidal constraints have
	loosely strong companions, then its underlying bicategory [of loose morphisms, N/A] $\clH(\dblD)$ is a monoidal bicategory. If $\dblD$ is braided or symmetric, so is $\clH(\dblD)$.
\end{proposition}
The verification that $\dblD$ is a monoidal double category consists of a lengthy computation, but a fairly more elementary one: `most of [the coherence diagrams] are fairly small [they are all 2-dimensional, A/N], and in any given case most or all of them are fairly obvious'.

The condition to check on the associator $\alpha,\lambda,\rho$ and unitors of the structure is then the following: [\emph{ibi}, 4.10]
\begin{itemize}
	\item $\alpha,\lambda,\rho$ all admit a companion (this will be implied by \autoref{has_cumpa} below);
	\item pasting the components of $\alpha,\lambda,\rho$ with the universal cells of their companions yields colax transformations that are in fact pseudo-natural.
\end{itemize}
These conditions are lengthy to write down but the most straightforward to check.
\subsection{Co/monads in \texorpdfstring{$\cartDblMly$}{Mly}}
\begin{definition}[Monad in a \PDC]\label{mnd_in_generic_pdc}
	Let $\dblD$ be a \PDC. A monad in $\dblD$ consists of a monad in the loose bicategory $\clH(\dblD)$:
	\begin{itemize}
		\item an object $A\in\dblD_0$ and a loose morphism $\fkm = (E,\langle d,s\rangle) : A\pto A$ (cf. \autoref{cart_dbl_mly});
		\item cells $\mu : \fkm\horC \fkm \To \fkm,\eta : \fki_A \To \fkm$ respectively called \emph{multiplication} and \emph{unit}, having shape
		      \[\begin{tikzcd}[cramped]
				      A & A & A & A & A \\
				      A && A & A & A
				      \arrow["\fkm", "\shortmid"{marking}, from=1-1, to=1-2]
				      \arrow[Rightarrow, no head, from=1-1, to=2-1]
				      \arrow["\fkm", "\shortmid"{marking}, from=1-2, to=1-3]
				      \arrow[Rightarrow, no head, from=1-3, to=2-3]
				      \arrow[""{name=0, anchor=center, inner sep=0}, "\fki_A", "\shortmid"{marking}, Rightarrow, no head, from=1-4, to=1-5]
				      \arrow[Rightarrow, no head, from=1-4, to=2-4]
				      \arrow[Rightarrow, no head, from=1-5, to=2-5]
				      \arrow[""{name=1, anchor=center, inner sep=0}, "\fkm"', "\shortmid"{marking}, from=2-1, to=2-3]
				      \arrow[""{name=2, anchor=center, inner sep=0}, "\fkm"', "\shortmid"{marking}, from=2-4, to=2-5]
				      \arrow["\mu"{description}, draw=none, from=1-2, to=1]
				      \arrow["\eta"{description}, draw=none, from=0, to=2]
			      \end{tikzcd}\]
	\end{itemize}
	and subject to the equations\footnote{Here and in many other instances, the diagrammatic conditions are zoomed out, using \autoref{def_pdc} as reference for the notation, and dropping the associator and unitor isomorphisms of the \PDC to save space. To save even more space, we will denote the conditions that cell compositions as $\frac{\mu\mid M}\mu = \frac{M\mid\mu}\mu$ or $\frac{M\mid\eta}\mu = M = \frac{\eta\mid M}\mu$ --the way in which these diagrams can be zoomed in being completely straightforward.}
	\kern-2em
\begin{tangle}{(8,3)}
	\tgBorderA{(0,0)}{white}{white}{white}{white}
	\tgBorder{(0,0)}{0}{1}{1}{0}
	\tgBorderA{(1,0)}{white}{white}{white}{white}
	\tgBorder{(1,0)}{0}{1}{0}{1}
	\tgBorderA{(2,0)}{white}{white}{white}{white}
	\tgBorder{(2,0)}{0}{1}{1}{1}
	\tgBorderA{(3,0)}{white}{white}{white}{white}
	\tgBorder{(3,0)}{0}{0}{1}{1}
	\tgBorderA{(4,0)}{white}{white}{white}{white}
	\tgBorder{(4,0)}{0}{1}{1}{0}
	\tgBorderA{(5,0)}{white}{white}{white}{white}
	\tgBorder{(5,0)}{0}{1}{1}{1}
	\tgBorderA{(6,0)}{white}{white}{white}{white}
	\tgBorder{(6,0)}{0}{1}{0}{1}
	\tgBorderA{(7,0)}{white}{white}{white}{white}
	\tgBorder{(7,0)}{0}{0}{1}{1}
	\tgBorderA{(0,1)}{white}{white}{white}{white}
	\tgBorder{(0,1)}{1}{1}{1}{0}
	\tgBorderA{(1,1)}{white}{white}{white}{white}
	\tgBorder{(1,1)}{0}{1}{0}{1}
	\tgBorderA{(2,1)}{white}{white}{white}{white}
	\tgBorder{(2,1)}{1}{1}{0}{1}
	\tgBorderA{(3,1)}{white}{white}{white}{white}
	\tgBorder{(3,1)}{1}{0}{1}{1}
	\tgBorderA{(4,1)}{white}{white}{white}{white}
	\tgBorder{(4,1)}{1}{1}{1}{0}
	\tgBorderA{(5,1)}{white}{white}{white}{white}
	\tgBorder{(5,1)}{1}{1}{0}{1}
	\tgBorderA{(6,1)}{white}{white}{white}{white}
	\tgBorder{(6,1)}{0}{1}{0}{1}
	\tgBorderA{(7,1)}{white}{white}{white}{white}
	\tgBorder{(7,1)}{1}{0}{1}{1}
	\tgBorderA{(0,2)}{white}{white}{white}{white}
	\tgBorder{(0,2)}{1}{1}{0}{0}
	\tgBorderA{(1,2)}{white}{white}{white}{white}
	\tgBorder{(1,2)}{0}{1}{0}{1}
	\tgBorderA{(2,2)}{white}{white}{white}{white}
	\tgBorder{(2,2)}{0}{1}{0}{1}
	\tgBorderA{(3,2)}{white}{white}{white}{white}
	\tgBorder{(3,2)}{1}{0}{0}{1}
	\tgBorderA{(4,2)}{white}{white}{white}{white}
	\tgBorder{(4,2)}{1}{1}{0}{0}
	\tgBorderA{(5,2)}{white}{white}{white}{white}
	\tgBorder{(5,2)}{0}{1}{0}{1}
	\tgBorderA{(6,2)}{white}{white}{white}{white}
	\tgBorder{(6,2)}{0}{1}{0}{1}
	\tgBorderA{(7,2)}{white}{white}{white}{white}
	\tgBorder{(7,2)}{1}{0}{0}{1}
	\tgCell{(1,0.5)}{\mu}
	\tgCell{(1.5,1.5)}{\mu}
	\tgCell{(2.5,0.5)}{\fkm}
	\tgCell{(3.5,1)}{=}
	\tgCell{(4.5,0.5)}{\fkm}
	\tgCell{(6,0.5)}{\mu}
	\tgCell{(5.5,1.5)}{\mu}
\end{tangle}\kern-1em and \kern-2em
\begin{tangle}{(9,3)}
	\tgBorderA{(0,0)}{white}{white}{white}{white}
	\tgBorder{(0,0)}{0}{1}{1}{0}
	\tgBorderA{(1,0)}{white}{white}{white}{white}
	\tgBorder{(1,0)}{0}{1}{1}{1}
	\tgBorderA{(2,0)}{white}{white}{white}{white}
	\tgBorder{(2,0)}{0}{0}{1}{1}
	\tgBorderA{(3,0)}{white}{white}{white}{white}
	\tgBorder{(3,0)}{0}{1}{1}{0}
	\tgBorderA{(4,0)}{white}{white}{white}{white}
	\tgBorder{(4,0)}{0}{1}{1}{1}
	\tgBorderA{(5,0)}{white}{white}{white}{white}
	\tgBorder{(5,0)}{0}{0}{1}{1}
	\tgBorderA{(6,0)}{white}{white}{white}{white}
	\tgBorder{(6,0)}{0}{1}{1}{0}
	\tgBorderA{(7,0)}{white}{white}{white}{white}
	\tgBorder{(7,0)}{0}{1}{0}{1}
	\tgBorderA{(8,0)}{white}{white}{white}{white}
	\tgBorder{(8,0)}{0}{0}{1}{1}
	\tgBorderA{(0,1)}{white}{white}{white}{white}
	\tgBorder{(0,1)}{1}{1}{1}{0}
	\tgBorderA{(1,1)}{white}{white}{white}{white}
	\tgBorder{(1,1)}{1}{1}{0}{1}
	\tgBorderA{(2,1)}{white}{white}{white}{white}
	\tgBorder{(2,1)}{1}{0}{1}{1}
	\tgBorderA{(3,1)}{white}{white}{white}{white}
	\tgBorder{(3,1)}{1}{1}{1}{0}
	\tgBorderA{(4,1)}{white}{white}{white}{white}
	\tgBorder{(4,1)}{1}{1}{0}{1}
	\tgBorderA{(5,1)}{white}{white}{white}{white}
	\tgBorder{(5,1)}{1}{0}{1}{1}
	\tgBorderA{(6,1)}{white}{white}{white}{white}
	\tgBorder{(6,1)}{1}{0}{1}{0}
	\tgBlank{(7,1)}{white}
	\tgBorderA{(8,1)}{white}{white}{white}{white}
	\tgBorder{(8,1)}{1}{0}{1}{0}
	\tgBorderA{(0,2)}{white}{white}{white}{white}
	\tgBorder{(0,2)}{1}{1}{0}{0}
	\tgBorderA{(1,2)}{white}{white}{white}{white}
	\tgBorder{(1,2)}{0}{1}{0}{1}
	\tgBorderA{(2,2)}{white}{white}{white}{white}
	\tgBorder{(2,2)}{1}{0}{0}{1}
	\tgBorderA{(3,2)}{white}{white}{white}{white}
	\tgBorder{(3,2)}{1}{1}{0}{0}
	\tgBorderA{(4,2)}{white}{white}{white}{white}
	\tgBorder{(4,2)}{0}{1}{0}{1}
	\tgBorderA{(5,2)}{white}{white}{white}{white}
	\tgBorder{(5,2)}{1}{0}{0}{1}
	\tgBorderA{(6,2)}{white}{white}{white}{white}
	\tgBorder{(6,2)}{1}{1}{0}{0}
	\tgBorderA{(7,2)}{white}{white}{white}{white}
	\tgBorder{(7,2)}{0}{1}{0}{1}
	\tgBorderA{(8,2)}{white}{white}{white}{white}
	\tgBorder{(8,2)}{1}{0}{0}{1}
	\tgCell{(7,1)}{\fkm}
	\tgCell{(0.5,0.5)}{\fkm}
	\tgCell{(4.5,0.5)}{\fkm}
	\tgCell{(1.5,0.5)}{\eta}
	\tgCell{(3.5,0.5)}{\eta}
	\tgCell{(1,1.5)}{\mu}
	\tgCell{(4,1.5)}{\mu}
	\tgCell{(2.5,1)}{=}
	\tgCell{(5.5,1)}{=}
\end{tangle}\kern-1em
	that exhibit, respectively, the \emph{associative} and \emph{unital} property for $\fkm$.
\end{definition}
When spelt out in $\cartDblMly$, the definition of a monad consists of the following data.
\begin{definition}[Monad in $\cartDblMly$]\label{deffa_monada}
	A monad in the \PDC of Mealy automata consists of a Mealy automaton
	$\footnotesize\begin{tikzcd}[cramped]
			E & {A\times E} & A
			\arrow["d"', from=1-2, to=1-1]
			\arrow["s", from=1-2, to=1-3]
		\end{tikzcd}$ arranged so that
	the unit cell $\eta : 1\to E$ and the multiplication cell $\mu$ respectively determine an element $e_0 \in E$ (as already pointed out in \autoref{gdaogao}) and a function $\mu : E\times E\to E$; these data are subject to the following requirements: first, the cell axioms for $\eta$ and $\mu$,
	\begin{enumerate}
		\item \label{ac_1} ($\eta\emdash s\emdash\mathrm{comp}$): $\forall a\in A : s (a , e_0) = a$;
		\item \label{ac_2} ($\eta\emdash d\emdash\mathrm{comp}$): $\forall a\in A : d(a,e_0) = e_0$.
		\item \label{mc_1} ($\mu\emdash s\emdash\mathrm{comp}$): $\forall a\in A \text{ and } e,e'\in E,\quad s (a , \mu(e , e')) = s (s (a , e) , e')$;
		\item \label{mc_2} ($\mu\emdash d\emdash\mathrm{comp}$): $\forall a\in A \text{ and } e,e'\in E,\quad d (a , \mu(e , e')) = \mu(d (a , e) , d (s (a , e) , e'))$.
	\end{enumerate}
	Furthermore, the triple $(E,e_0,\mu)$ forms a monoid --the multiplication of which we will from now on denote as an infix dot, $e\cdot e'=\mu(e,e')$:
	\begin{enumtag}{ma}
		\item\label{ma_1} $\mu$ is associative, i.e. for all $e,e',e''\in E$, $(e\cdot e')\cdot e'' = e\cdot (e'\cdot e'')$;
		\item\label{ma_2} $\mu$ has $e_0$ as neutral element, $e\cdot e_0 = e = e_0\cdot e$.
	\end{enumtag}
\end{definition}
\begin{remark}[Canonical extensions of $d,s$]\label{s_extension}
	In the same notation as above, denote $d^+(as,e)$ the \emph{canonical extension} of $d : A\times E \to E$ to an action of the monoid $A^*$ of lists on $A$ over the set $E$, defined inductively as $d^+(\emptyList,e):=e$, $d^+(a\cons as,e):=d(a, d^+(as, e))$.

	In the following, we denote $d^+,s$ as infix symbols
	$d^+(as,e)=as\dact e$ and $s(a,e)=a\odot e$, so that condition \ref{mc_2} rewrites as
	\[\label{fuffa_magica}a\dact (e\cdot e') = (a \dact e) \cdot ((a\odot e) \dact e')\]
	Note that there is also a canonical extension of $s$ to a map $A^*\times E\to A^*:(as,e)\mapsto as\sact e$ defined inductively as $\emptyList \sact e = \emptyList$ and $(a \cons as) \sact e = s (a , as \dact e) \cons (as \sact e)$ (so $(a\cons\emptyList)\sact e = s(a,e)\cons\emptyList$ by definition, and $\sact$ really extends $\odot$); moreover, \eqref{fuffa_magica} holds more in general for $\sact$ in the form
	\[\label{fuffa_magica2}as\dact (e\cdot e') = (as \dact e) \cdot ((as\sact e) \dact e').\]
\end{remark}
This provides the following characterization of a monad in $\cartDblMly$ in terms of two `matched' monoid actions
\begin{corollary}[Monads as matching pairs]\label{sgrulbio}
	A monad $(\fkm,\mu,\eta)$ in $\cartDblMly$ as in \autoref{deffa_monada}, on the object $A$, consists of a pair of monoid actions:
	\begin{itemize}
		\item an action $d^* : A^*\times E\to E$, obtained as the canonical extension of $d$;
		\item an action of $E$ (which is a monoid due to \ref{ma_1}, \ref{ma_2}) on the set $A$, obtained as the map $s : A\times E\to A$ (\autoref{deffa_monada}.(\ref{ac_1}), \autoref{deffa_monada}.(\ref{mc_1}) say that it is an action).
	\end{itemize}
	These actions are subject to the following mutual compatibility conditions:
	\begin{itemize}
		\item the identity element $e_0$ in the monoid $E$ is a fixpoint for the action $d^*$;
		\item dually, the identity element $\emptyList$ in $A^*$ is a fixpoint for $\sact$;
		\item the two \emph{bicrossed equations} \eqref{fuffa_magica2} (relating the action $d^*$ and the monoid structure of $E$) and
		      \[\label{altra_fuffa}(as \concat bs) \sact e = (as \sact (bs \dact e)) \concat (bs \sact e)\]
		      (relating the action $\sact$ and the monoid structure of $A^*$) hold.\footnote{Some authors, especially when defining the analogue notion in group theory, refer to two monoids $G,H$, equipped with the actions compatible in that way, as being \emph{matched}, or forming a \emph{matching pair} $(G,H;\dact,\sact)$.}
	\end{itemize}
\end{corollary}
(\autoref{altra_fuffa} can be proved by induction.) It turns out that the data specified by \autoref{sgrulbio} is well-known to algebraists: cf. \cite{Agore_2014,Agore2009forfinite,Agore_2009,Agore_2010}. Whenever two monoids $G,H$ act one over the other with maps satisfying the bicrossed equations, the \emph{bicrossed product} $H\bowtie G$ can be constructed, generalizing their semidirect product (in the sense that the latter is obtained as particular instance of the former).\footnote{In both the group and monoid case, when $\sact$ (resp., $\dact$) is the trivial `constant at 1' action, then $\dact$ (resp., $\sact$) is an action via endomorphisms (or group automorphisms) of $G$ on $H$ (resp., of $H$ on $G$), and $H\bowtie G$ is the semidirect product $H\rtimes G$(resp., $H\ltimes G$).}
\begin{definition}[Bicrossed product of monoids]\label{crucial_constru}
	Let $G,H$ be two monoids (whose multiplication we denote conveniently as $\concat$ and $\cdot$ respectively, as $G=E$ and $H=A^*$ in the case we're interested), and consider
	\begin{itemize}
		\item a left action of $G$ on $H$, denoted $\firstblank\dact\firstblank : G\times H \to H$;
		\item a right action of $H$ on $G$, denoted $\firstblank\sact\firstblank : G\times H \to G$;
		\item such that $g\dact h = 1_H$ and $g\sact h = 1_G$, and
		\item satisfying the \emph{bicrossed equations} \eqref{fuffa_magica2} and \eqref{altra_fuffa},
		      \[\notag as\dact (e\cdot e') = (as \dact e) \cdot ((as\sact e) \dact e')\qquad (as \concat bs) \sact e = (as \sact (bs \dact e)) \concat (bs \sact e)\]
	\end{itemize}
	Then the \emph{bicrossed product} $H\bowtie G$ of the two monoids is defined as the monoid having
	\begin{itemize}
		\item as carrier the Cartesian product $H\times G$;
		\item identity element the pair $(1_H,1_G)$;
		\item monoid operation $\bullet$ defined using \emph{both} actions as
		      \[
			      \label{mult_in_bicro}
			      (x , as) \bullet (y , bs) =
			      \big(x \cdot (as \dact y) , (as \sact y) \concat bs \big)
		      \]
	\end{itemize}
\end{definition}
The proof is omitted, as it is easy to adapt from \cite[Theorem 1.3]{Agore2009forfinite} where it is stated for groups; unitality is straightforward, and associativity crucially depends on the bicrossed equations. The only non-immediate part proving the following statement is to establish that \eqref{altra_fuffa} holds for the action $\sact$ of $E$ on $A^*$, but this is an easy exercise on induction for lists.
\begin{corollary}\label{dgaddaaa}
	Every monad in $\cartDblMly$ packages the structure necessary to define the `semifree' bicrossed product $E\bowtie A^*$ between the free monoid on $A$ and the monoid $(E,\cdot,e_0)$ obtained from \ref{ma_1}, \ref{ma_2}.\footnote{This is actually part of a functor $\beta : \Mnd(\cartDblMly) \to \Mon$ sending $(A\overset M\pto A,\mu,\eta)\mapsto E\bowtie_M A^*$, for the natural notion of monad morphism $(M,\mu,\eta) \to (N,\nu,\gamma)$, that we recall in \autoref{mormonad_tig}.}
\end{corollary}
The paper \cite[Proposition 1.7.(i)]{Agore2009forfinite} states a colimit-like universal property for the bicrossed product of groups; a similar one for monoids is easily stated, and the proof is exactly as in \cite[Proposition 1.7.(i)]{Agore2009forfinite}.
\begin{theorem}\label{sbriggosbraggo}
	The monoid $H\bowtie G$, equipped with the obvious injections $H \xto{i_H}H\bowtie G\xot{i_G} G$ is the initial object of the category where an object is a cospan $H \xto{u_H}N\xot{u_G} G$ of monoid homomorphisms such that $u_G(g)\bullet u_H(h) = (g\dact h,g\sact h)$ and $u_H(h)\bullet u_G(g)=(h,g)$.
\end{theorem}
It is easy to deduce from this universal property a characterization of $H\bowtie G$ as a coequalizer, the cokernel of a certain map $H*G\to H*G$ induced via $i_H,i_G$; however, such presentation for $H\bowtie G$ is not particularly enlightening. Instead, it's easy to characterise the presheaf category $\Set^{H\bowtie G}$ in terms of the categories $\Set^H,\Set^G$.
\begin{lemma}\label{crucialemma}
	(Proof below) Let $(G,H;\dact,\sact)$ be a matching pair as in \autoref{crucial_constru}; then it is equivalent to give
	\begin{enumerate}
		\item a (left) representation $(X,\star)$ of the bicrossed product $H\bowtie G$, i.e. a functor $H\bowtie G \to \Set$ inducing $\firstblank\star\firstblank : (H\bowtie G)\times X\to X$;
		\item a pair of (left) representations $\alpha : H \times X \to X$, $\beta : G\times X \to X$ subject to the equation \[\label{gndgoadbgoadb}\forall x\in X : \beta(g, \alpha(h,x)) = \alpha(g\dact h, \beta(g\sact h,x))\]
	\end{enumerate}
	Stated in string diagrams, the condition \eqref{gndgoadbgoadb} reads
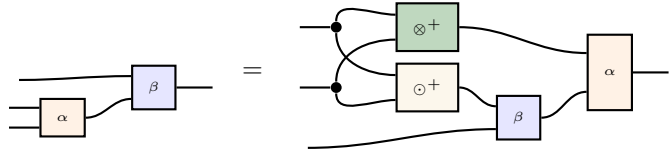
	\begin{center}
		\begin{tikzpicture}[scale=0.4]
			\node[font=\tiny,box=0/1/0/2,fill=orange!10] (alpha) at (0, -1) {\alpha};
			\node[font=\tiny,box=0/1/0/2,fill=blue!10] (beta) at (3, 0) {\beta};
			\node (w) at (5.25, 0) {};
			\node (hanging) at (-1.75,0.25) {};
			\wires[thick]{
				alpha = { east = beta.west.2 },
				beta = { east = w.west },
				hanging = { east = beta.west.1 }
			}{ alpha.west.1 , alpha.west.2 }
			\begin{scope}[xshift=9cm]
				\node[dot] (x) at (0, 2) {};
				\node[dot] (y) at (0, 0) {};
				\node (z) at (-1.25, -2) {};
				\node (w) at (11.25, .5) {};
				\node[font=\tiny,box=0/1/0/2,fill=dColor!25] (dact) at (3,2) {\dact};
				\node[font=\tiny,box=0/1/0/2,fill=sColor!25] (sact) at (3,0) {\sact};
				\node[font=\tiny,box=0/1/0/2,fill=blue!10] (beta') at (6,-1) {\beta};
				\node[font=\tiny,box=0/1/0/2,minimum height=1cm,fill=orange!10] (alpha') at (9,.5) {\alpha};
				\wires[thick]{
				x = { north = dact.west.1 , south = sact.west.1 },
				y = { north = dact.west.2 , south = sact.west.2 },
				dact = { east = alpha'.west.1 },
				sact = { east = beta'.west.1 },
				beta' = { west.2 = z.east , east = alpha'.west.2 },
				alpha' = { east = w.west },
				}{ x.west , y.west }
			\end{scope}
			\node at (6.25,.5) {$=$};
		\end{tikzpicture}
	\end{center}
	Given two actions $\alpha,\beta$ as above, posing $(h,g)\star x := \alpha(h, \beta(g,x))$ defines an action of $H\bowtie G$.
\end{lemma}
\subsection{Modules over a monad}
\autoref{dgaddaaa} suggests the natural conjecture: algebras for a monad in $\cartDblMly$ (and thus, in $\proc\clK$) coincide with sets equipped with an action of the bicrossed product. This is indeed the case.
\begin{definition}[Left modules for a monad on a \PDC]
	Let $(M : A\pto A,\mu,\eta)$ be a monad on an object $A$ of a \PDC $\dblD$ as in \autoref{sgrulbio}; a (left) \emph{$M$-module}, for $M$ consists of a pair $(P,\xi)$ where $P : A\pto X$ is a loose morphism, and $\xi$ is a cell
	\[\begin{tikzcd}[cramped]
			A & A & X \\
			A && X
			\arrow["M", "\shortmid"{marking}, from=1-1, to=1-2]
			\arrow[Rightarrow, no head, from=1-1, to=2-1]
			\arrow["\xi"{description}, draw=none, from=1-1, to=2-3]
			\arrow["P", "\shortmid"{marking}, from=1-2, to=1-3]
			\arrow[Rightarrow, no head, from=1-3, to=2-3]
			\arrow["P"', "\shortmid"{marking}, from=2-1, to=2-3]
		\end{tikzcd}\]
	subject to the two axioms of unitality and associativity: $\frac{\eta \mid P}{\xi} = P$ and $\frac{M\mid\xi}{\xi} = \frac{\mu\mid P}{\xi}$.
\end{definition}
\begin{theorem}[Characterization of left modules for monads in $\cartDblMly$]\label{crucial_thm} (Proof below)
	There is an isomorphism of categories between
	\begin{itemize}
		\item the category of left modules for a monad $(\fkm,\mu,\eta)$ on $A$;
		\item the category of $\Set$-actions of the monoid $E\bowtie_M A^*$, obtained in \autoref{dgaddaaa}.
	\end{itemize}
\end{theorem}
\begin{remark}\label{bonbiollo}
	If $\clK$ is Cartesian, the theory of \emph{comonads} is instead rather trivial, as a consequence of Fox theorem, asserting that every object of a Cartesian category is equipped with a canonical choice of comonoid, for which every morphism is a homomorphism. The notion of comonoid (and Hopf monoid of some sort, assumption that becomes necessary when restating the matching condition) can evidently become very interesting in nonCartesian monoidal categories (cf., for example, \cite{loregian2024automata}, on the lines of \cite{adam-trnk:automata}, for a nonCartesian category with an interesting -generalized- automata theory).

	The matching conditions have an evident interpretation in monoidal categories where the following string diagrams live:
	
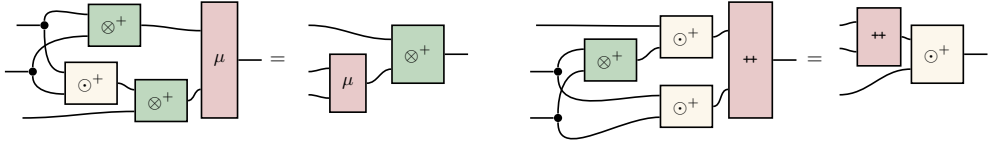
\begin{adju}
		\begin{tikzpicture}[scale=0.4]
			\node[dot] (x) at (0, 2) {};
			\node[dot] (y) at (-.5, 0) {};
			\node (z) at (-1.25, -2) {};
			\begin{scope}[xshift=2cm]
				\node[font=\small,box=0/1/0/2,fill=dColor!25] (otimes1) at (1,2) {\dact};
				\node[font=\small,box=0/1/0/2,fill=sColor!25] (odot) at (0,-.5) {\sact};
				\node[font=\small,box=0/1/0/2,fill=dColor!25] (otimes2) at (3,-1.25) {\dact};
				\node[font=\small,box=0/1/0/2, minimum height=2cm,fill=muColor!25] (mu) at (5.5,.5) {\mu};
			\end{scope}
			\wires[thick]{
				x = {north = otimes1.west.1, south = odot.west.1 },
				y = {north = otimes1.west.2, south = odot.west.2 },
				odot = { east = otimes2.west.1 },
				otimes1 = { east = mu.west.1 },
				otimes2 = { east = mu.west.2, west.2 = z.east },
			}{ x.west, y.west, mu.east }
			\begin{scope}[xshift=11cm]
				\node (x2) at (0, 2) {};
				\node (y2) at (0, 0) {};
				\node (z2) at (0, -1) {};
				\node[font=\small,box=0/1/0/2,minimum height=1cm,fill=muColor!25] (mu2) at (2,-.5) {\mu};
				\node[font=\small,box=0/1/0/2,minimum height=1cm, fill=dColor!25] (otimes3) at (5,.75) {\dact};
				\wires[thick]{
					mu2 = { west.2 = z2.east , west.1 = y2.east , east = otimes3.west.2 },
					otimes3 = { west.1 = x2.east }
				}{ otimes3.east }
			\end{scope}
			\node at (10,.5) {$=$};
			\begin{scope}[xshift=22cm]
				\node (x) at (-1.25, 2) {};
				\node[dot] (y) at (0, 0) {};
				\node[dot] (z) at (0, -2) {};
				\begin{scope}[xshift=2cm]
					\node[font=\small,box=0/1/0/2,fill=sColor!25] (odot1) at (3.5,1.5) {\sact};
					\node[font=\small,box=0/1/0/2,fill=dColor!25] (otime) at (.25,.5) {\dact};
					\node[font=\small,box=0/1/0/2,fill=sColor!25] (odot2) at (3.5,-1.5) {\sact};
					\node[font=\small,box=0/1/0/2, minimum height=2cm,fill=muColor!25] (concat) at (6.25,.5) {\concat};
				\end{scope}
				\wires[thick]{
					odot1 = { west.1 = x.east , east = concat.west.1 },
					odot2 = { east = concat.west.2 },
					otime = { east = odot1.west.2 },
					y = { north = otime.west.1 , south = odot2.west.1 },
					z = { north = otime.west.2 , south = odot2.west.2 },
				}{ y.west , z.west , concat.east }
				\begin{scope}[xshift=11.75cm]
					\node (z2) at (0, 2) {};
					\node (y2) at (0, 1) {};
					\node (x2) at (0, -1) {};
					\node[font=\small,box=0/1/0/2,minimum height=1cm,fill=muColor!25] (concat2) at (2,1.5) {\concat};
					\node[font=\small,box=0/1/0/2,minimum height=1cm, fill=sColor!25] (otimes3) at (4.5,.75) {\sact};
					\wires[thick]{
						concat2 = { west.1 = z2.east , west.2 = y2.east , east = otimes3.west.1 },
						otimes3 = { west.2 = x2.east }
					}{ otimes3.east }
				\end{scope}
				\node at (11,.5) {$=$};
			\end{scope}
		\end{tikzpicture}
	\end{adju}
	It becomes of great interest, for example, to form the bicrossed product $H\bowtie H^\op$ on a Hopf algebra $H$ with its opposite, \cite{Drinfeld1988,Majid1990,Muger_2008}; the same construction was generalized to categories \cite{Kassel_1995,Majid1991} --a direction that doesn't feel peregrine to explore in light of the results in \cite{guitart1974remarques,guitart1978bimodules,Guitart1977} where a bicategory of Mealy automata in the ambient category $\Cat$ is introduced; Majid in \cite{Majid1990,Majid1991} studies bicrossed products of monoidal categories, and thus constitutes a vertical categorification of the theory outlined here; cf. also \cite{Street_1998}.
\end{remark}
\begin{remark}[Functoriality in \autoref{dgaddaaa}]
	The above characterization of monads in $\cartDblMly$ is part of an equivalence of categories:
	\begin{itemize}
		\item the category $\cate{DblMnd}(\cartDblMly)$ ($\cate{DblMnd}$ for short) of double monads has objects the monads $(\fkm=(E,\langle d,s\rangle) : A \pto A,\eta,\mu)$ and morphisms of monads $(u,f) : (\fkm,\eta_A,\mu_A) \to (\fkn,\eta_B,\mu_B)$ those of \autoref{mormonad_tig};
		\item the category $\cate{XXPair}$ of bicrossed pairs has objects the bicrossed pairs $(E,A,\dact,\sact)$ as in \autoref{crucial_constru} and morphisms $(E,A,\dact_1,\sact_1) \to (E',B,\dact_2,\sact_2)$ those pairs of monoid homomorphisms $f : E\to E'$, $u : A\to B$ suitably compatible with actions, in the sense that the diagram
		      $\footnotesize\begin{tikzcd}[cramped]
				      E & {A\times E} & A \\
				      {E'} & {B\times E'} & B
				      \arrow["f"', from=1-1, to=2-1]
				      \arrow["{\dact_1}"', from=1-2, to=1-1]
				      \arrow["{\sact_1}", from=1-2, to=1-3]
				      \arrow["{u\times f}"', from=1-2, to=2-2]
				      \arrow["u", from=1-3, to=2-3]
				      \arrow["{\dact_2}", from=2-2, to=2-1]
				      \arrow["{\sact_2}"', from=2-2, to=2-3]
			      \end{tikzcd}$
		      commutes. This is enough to induce a monoid homomorphism $u\bowtie f : E\bowtie A \to E' \bowtie B'$.
	\end{itemize}
\end{remark}
The proof of the following statement is straightforward:
\begin{proposition}[Monads are semifree bicrossed products]
	There is a functor $\cate{DblMnd} \to \cate{XXPair}$ sending a monad $(\fkm,\eta,\mu)$ to the matching pair determined as in \autoref{sgrulbio} (and thus to the bicrossed product of the monoids $A^*$ and $E$); a monad morphism corresponds precisely to a morphism in $\cate{XXPair}$ where the homomorphism $u : A^* \to B^*$ is actually induced by a function on the generators $u_0 : A \to B$; this identifies $\cate{DblMnd}$ with the nonfull subcategory $\cate{SFXPair}$ of $\cate{XXPair}$ spanned by matching pairs whose second component is free, and homomorphisms on that component are induced on the generators.
\end{proposition}
Another consequence of \autoref{sgrulbio} is that it makes transparent the construction of the free monad on a loose endomorphism $F : A \pto A$:
\begin{corollary}[Free monad on a loose endomorphism]\label{coro_fremonad}
	From \autoref{sgrulbio} we deduce that a loose endomorphism $F : A\pto A$ in $\cartDblMly$ is given by a `matching pair of sets'
	\[\begin{tikzcd}[cramped]
			E & {A\times E} & A
			\arrow["d"', from=1-2, to=1-1]
			\arrow["s", from=1-2, to=1-3]
		\end{tikzcd}\]
	The functions $d,s$ have a canonical extension $M(F)$
	\[\begin{tikzcd}[cramped]
			E^* & {A\times E^*} & A
			\arrow["d^+"', from=1-2, to=1-1]
			\arrow["s^+", from=1-2, to=1-3]
		\end{tikzcd}\]
	defined as follows:
	\[\begin{cases}
			d^+(a , \emptyList) = \emptyList \\
			d^+(a , e \cons es) = d (a , e) \cons d^+ (a , es)
		\end{cases}\qquad
		\begin{cases}
			s^+(a , \emptyList) = a \\
			s^+(a , e \cons es) = s^+ (s (a , e) , es^\leftarrow)
		\end{cases}\]
	where $es^\leftarrow$ is the \emph{safe reverse} of a list --empty if $es=\emptyList$.
	One can easily check that this extension $M(F)$ provides the \emph{free monad} on the loose morphism $F$, i.e.:
	\begin{itemize}
		\item $M(F)$ is equipped with cells $\eta,\mu$ as in the definition of monad in \ref{deffa_monada}, defining a monad in the \PDC of Mealy automata, and with a cell $\smallCell F{}{}{M(F)}[\nu]$;
		\item the universal property of a free object holds: given a cell $\smallCell F{}{}N[\gamma]$ where $N$ is a double monad, there is a unique monad morphism (cf. \autoref{mormonad_tig}) $\smallCell {M(F)}{}{}N[\gamma^*]$ such that $\gamma^*\vertC \nu = \gamma$.
	\end{itemize}
\end{corollary}
This entire construction assembles into a double functor, providing a left adjoint to the forgetful, sending a double monad to its underlying loose endomorphism $M$.
\begin{remark}[There are more monads, Horatio...]\label{horatio}
	Another question is: how do monads relate to adjunctions in the \PDC of Mealy automata? The answer is, relatively poorly: the notion of loose adjunction (cf. \autoref{def_loose_adja}) entails that in order for a pair of loose cells to fit into an adjunction $\ell\dashv r$, i.e. to be equipped with unit and counit maps satisfying the zig-zag identities, the carrier of both automata
	\[\begin{tikzcd}[cramped]
			E & {A\times E} & B & {E'} & {B\times E'} & A
			\arrow["{d_\ell}"', from=1-2, to=1-1]
			\arrow["{s_\ell}", from=1-2, to=1-3]
			\arrow["{d_r}"', from=1-5, to=1-4]
			\arrow["{s_r}", from=1-5, to=1-6]
		\end{tikzcd}\]
	must be singletons; thus, the left adjoint $\ell : A \pto B$ must be the companion $(f_\ell)_*$ of a tight morphism $f_\ell : A\to B$. This, coupled with the fact that then $r = f_\ell^*$ is the conjoint of $f$, means that the unique adjunctions are induced as companion-conjoint pairs $f_*\dashv f^*$ (cf. \cite[4.1.4]{Grandis2019} for a proof that the following are equivalent: (i) $f : A\to B$ admits a companion $f_*$ and a conjoint $f^*$; (ii) $f_*\dashv f^*$; (iii) $f$ has a companion $f_*$ and $f_*$ is a left adjoint). Each such adjunction induces a monad for trivial reasons, $M : A\times 1\to 1\times A$ (and clearly, $1\bowtie A^*\cong A^*$); there are, however, plenty of other monads!
\end{remark}
We conclude the section discussing the structure of loose \emph{monad maps}: unravelling the definition relates to \cite[§2]{guitart1974remarques} and to our \cite[2.12---2.18]{EPTCS397.1}, and give some more motivation for the additional condition postulated in \cite{guitart1974remarques} to embed Mealy automata between monoids in a bicategory of spans.
\subsection{Double monad maps}\label{dbl_moma}
\begin{definition}[Loose monad map in a \PDC]
	Let $\dblD$ be a \PDC, and $M : A\pto A,N : B\pto B$ two monads as in \autoref{deffa_monada}; a \emph{monad map} $(U,\lambda) : M\to N$ consists of a loose morphism $U : A\pto B$, and of a cell
	\[\label{dist_diag}\begin{tikzcd}[cramped]
			A & A & B \\
			A & B & B
			\arrow["M", "\shortmid"{marking}, from=1-1, to=1-2]
			\arrow[Rightarrow, no head, from=1-1, to=2-1]
			\arrow["\lambda"{description}, draw=none, from=1-1, to=2-3]
			\arrow["U", "\shortmid"{marking}, from=1-2, to=1-3]
			\arrow[Rightarrow, no head, from=1-3, to=2-3]
			\arrow["U"', "\shortmid"{marking}, from=2-1, to=2-2]
			\arrow["N"', "\shortmid"{marking}, from=2-2, to=2-3]
		\end{tikzcd}\]
	subject to the following compatibility with the units $\eta^M,\eta^N$ and multiplications $\mu^M,\mu^N$ of $M,N$ (unitors and associators have been suppressed): $\frac{\mu^M\mid U}\lambda = \frac{\textstyle\frac{M\mid \lambda}{\lambda\mid N}}{U\mid \mu^N} $ and $\frac{\eta^M\mid U}{\lambda} = U\mid \eta^N$.
\end{definition}
\begin{definition}[Loose monad map in $\cartDblMly$]\label{mormonad_loo}
	Unwinding the definition, a monad map from $\fkm : A\times E \mathrel{\overset{\langle d_M,s_M\rangle}\to} E\times A$ to $\fkn : B\times E\mathrel{'\overset{\langle d_N,s_N\rangle}\to} E'\times B$
	in $\cartDblMly$ amounts to
	\begin{itemize}
		\item A Mealy automaton $\fku : A\times X \mathrel{\overset{\langle d_U,s_U\rangle}\to} X\times B$ fitting in a diagram like \eqref{dist_diag},
		\item whose structure cell $\lambda$ defines in turn a function $\langle \delta,\sigma\rangle : E \times X \to X \times E'$, i.e. a new Mealy automaton $\fkw$ between the state spaces of $M$ and $N$.
	\end{itemize}
	The automaton's structure maps are subject to the following conditions of compatibility with the structure maps of $M$ and $N$:
	\begin{enumtag}{dl}
		\item \label{dl_1} cell equations for $\lambda$:
		\begin{gather*}
			s_N (s_U (a , \delta (e , x)) , \sigma (e , x)) = s_U (s_M (a , e) , x) \\
			\begin{cases}
				\delta (d_M (a , e) , d_U (s_M (a , e) , x)) = d_U (a , \delta (e , x)) \\
				\sigma (d_M (a , e) , d_U (s_M (a , e) , x)) = d_N (s_U (a , \delta (e , x)) , \sigma (e , x))
			\end{cases}
		\end{gather*}
		\item \label{dl_2} compatibility with the units of $M,N$: a condition that splits into the two
		\[\delta (e_0 , x) = x \qquad \qquad \sigma (e_0 , x) = e_0'\]
		\item \label{dl_3} compatibility with the multiplications of $M,N$: a condition that splits into the two
		\[\oldstylenums{1}) : \delta (\mu_M (e , e') , x) = \delta (e , \delta (e' , x))\quad
			\oldstylenums{2}) : \sigma (\mu_M (e , e') , x) = \mu_N (\sigma (e , \delta (e' , x)) , \sigma (e' , x)).\]
	\end{enumtag}
\end{definition}
The second condition in \ref{dl_3} asserts a compatibility between the output function $\sigma$ of $\fkw$ and the monoid structure of $\fkm$'s and $\fkn$'s state spaces; natural conditions on Mealy automata between monoids are not new: \cite[§2]{guitart1974remarques} studies automata between monoids whose $s$ leg satisfies precisely condition \ref{dl_3}.\oldstylenums{2},\footnote{Guitart does not give a name to \ref{dl_3}.\oldstylenums{2}; in \cite{EPTCS397.1} we call it \emph{fugality} for $s$ (cf. [\emph{ibi}] for a thorough description of the matter); the condition can be motivated with the fact that fugality for $s$ is the only condition needed to induce a functor $\clE[d^+]\to B$ through $s$, where $\clE[d^+]$ is the domain of the discrete opfibration over $A^*$ naturally associated to the representation $d^+$.} because the bicategory they form fully embeds into a bicategory of spans; Guitart then proceeds to characterize the latter category as the Kleisli bicategory of a certain KZ monad of free cocompletion. In \cite{EPTCS397.1} we refurbish and extend such a result proving that \cite{guitart1974remarques}'s result defines a biadjunction [\emph{ibi}, Theorem 2.18] and that the sub-bicategory of fugal Mealy automata 1-fully and 2-fully embeds into Guitart's bicategory of spans, thus really providing a span representation for Mealy automata.
\begin{remark}[A double category of monads]
	From the notion of a monad map, we infer the notion of distributive law; these turn out to be interesting but unwieldy to spell out explicitly.

	We summarize here a point for future investigation: the tight monads maps of \autoref{mormonad_tig} and the loose ones of \autoref{mormonad_loo} above can be organized into a double category, for a suitable choice of cells with boundary
	$\footnotesize\begin{tikzcd}[cramped]
			M & {M'} \\
			N & {N'}
			\arrow["U", "\shortmid"{marking}, from=1-1, to=1-2]
			\arrow["f"', from=1-1, to=2-1]
			\arrow["\varpi"{description}, draw=none, from=1-1, to=2-2]
			\arrow["g", from=1-2, to=2-2]
			\arrow["V"', "\shortmid"{marking}, from=2-1, to=2-2]
		\end{tikzcd}$,
	whence one can find that a distributive law between monads is a (loose) monad in this bicategory of monads, $N : M\pto M$.
\end{remark}

\begin{remark}\label{how_translate}
	Applying \autoref{shllo_rula} does not make use of Cartesianity; a monad in $\dblMly$ consists now of a monoid $(E, e_0, \mu)$ internal to $(\clK,\otimes,I)$, equipped with a map $m : A\otimes E \to E\otimes A$, related to $\mu$ by the cell equations \autoref{baiadiai} (now keeping in mind that a 1-cell does not split as a span like in \autoref{florpo}):
	\begin{center}
		\begin{tikzpicture}[scale=0.4]
			\node (l) at (-2,-2) {};
			\node (r) at (8,-1.5) {};
			\node[inner sep=0,font=\tiny,box=0/2/0/2,fill=sColor] (m1) {m};
			\node[inner sep=0,font=\tiny,box=0/2/0/2,fill=sColor] at (3,-2) (m2) {m};
			\node[inner sep=0,font=\tiny,box=0/1/0/2,fill=red!20] at (6,-.5) (mu) {\mu};
			\wires[thick]{
				m1 = { east.1 = mu.west.1 , east.2 = m2.west.1 },
				m2 = { east.1 = mu.west.2 , west.2 = l.east , east.2 = r.west },
			}{ m1.west.1 , m1.west.2 , mu.east }
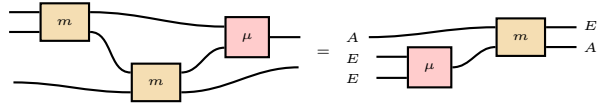
			\begin{scope}[xshift=9cm, yshift=.5cm]
				\node[left, font=\tiny] (l') at (1,-1) {$A$};
				\node[below=.5em of l', font=\tiny] {$E$};
				\node[below=-.25em of l', font=\tiny] {$E$};
				\node[font=\tiny, right] (E) at (7.75,-.625) {$E$};
				\node[font=\tiny, below=-.25em of E] {$A$};
				\node[font=\tiny,box=0/1/0/2,fill=red!20] at (3,-2) (mu') {\mu};
				\node[font=\tiny,box=0/2/0/2,fill=sColor] at (6,-1) (m') {m};
				\wires[thick]{
				mu' = { east = m'.west.2 },
				m' = { west.1 = l'.east }
				}{ mu'.west.1 , mu'.west.2 , m'.east.1 , m'.east.2 }
			\end{scope}
			\node[font=\tiny] at (8.5,-1) {$=$};
		\end{tikzpicture}
	\end{center}
	A monad in this context looks like the lax version (?) of a Drinfeld centralizer's component \cite{Joyal1991yang,Muger_2008}.

	Characterizing monads and their modules in this context (and the double category of monads \autoref{dbl_moma}) requires diving deep into representation theoretic phenomena.
\end{remark}
\section{Tight co/limits in \texorpdfstring{$\cartDblMly$}{Mly}}
The theory of co/limits in a \PDC is outlined in \cite{CTGDC_1999__40_3_162_0}, which we employ as reference for the notation; see [\emph{ibi}, §5] and in particular [\emph{ibi}, Theorem 5.5] asserting that small limits can be constructed from (tightly discrete) small products, (tightly discrete) equalizers and tabulators. We will leverage this result to prove that
\begin{theorem}[Existence of tight limits]
	The \PDC $\cartDblMly$ has products and equalizers -hence, it has all tightly discrete limits, provided $\clK$ admits limits; it does not admit tabulators.
\end{theorem}
and dually, that
\begin{theorem}[Existence of tight limits]
	The \PDC $\cartDblMly$ does not have coproducts or coequalizers -hence, it lacks all tightly discrete colimits, no matter their presence in $\clK$; it admits cotabulators.
\end{theorem}
We will split the argument into various steps. The simplest line of reasoning is that
\begin{remark}\label{coprod_no}
	Initial objects and coproducts do not exist in $\cartDblMly$.
\end{remark}
\begin{savvyProof}
	The non-existence of an initial object follows from \autoref{a_trick}; such an object corresponds to a unique cell $u_F$ so that the square
	$\footnotesize\begin{tikzcd}[cramped, ampersand replacement=\&]
			{\emptyset\times 1} \& {1\times \emptyset} \\
			{A\times E} \& {E\times B}
			\arrow[from=1-1, to=1-2]
			\arrow["{!\times e_0}"', from=1-1, to=2-1]
			\arrow["{e_0\times !}", from=1-2, to=2-2]
			\arrow[from=2-1, to=2-2]
		\end{tikzcd}$
	commutes; this forces $\emptyset$ to be an initial object in $\clK$, but also forces $E$ to be a singleton. For a similar reason, just a little bit more convoluted to enforce, coproducts, in the sense of \autoref{def_procopro_pc} (see also \cite{CTGDC_1999__40_3_162_0}) do not exist.
\end{savvyProof}
Moreover, the universal property of terminal objects and pullbacks amounts exactly to the double categorical universal property in $\cartDblMly$, thus:
\begin{theorem}\label{lims_yes}
	$\cartDblMly$ admits a (tight) terminal object and (tight) pullbacks, in the sense of \autoref{def_initerm_pc}, \autoref{def_pullpush_pc}, as soon as $\clK$ does admit a terminal object and pullbacks. (The same line of reasoning proves it admits wide pullbacks, or products of all cardinalities and equalizers if $\clK$ does.)
\end{theorem}
\subsection{Some genuinely double co/limits}\label{genuinely}
A direct computation proves that $\cartDblMly$ has all cotabulators, but almost no tabulators at all: so, automata have a \emph{collage} construction (`heteromorphisms are representable', cf. the example of enriched categories in \cite[3.7.1., 4.5.5.]{Grandis2019}) but not a \emph{category of elements} construction.
\begin{example}\label{constr_cotab}
	The cotabulator of a Mealy automaton $\fkf = \langle d,s\rangle : A\times E \to E\times B$ in $\cartDblMly$ is defined as follows:
	\begin{itemize}
		\item the carrier $\cotabu \fkf$ is defined from the coproduct $A+B$ as the coequalizer of  the two parallel arrows in
		      \[\begin{tikzcd}[cramped]
				      {A\times E} && {A+B} & {\cotabu \fkf}
				      \arrow["{\inj_B\circ s}", shift left=1, from=1-1, to=1-3]
				      \arrow["{\inj_A\circ \pi_A}"', shift right=1, from=1-1, to=1-3]
				      \arrow[from=1-3, to=1-4, "q"]
			      \end{tikzcd}\]
		      it is, in fact, easy to check that every cell like $\xi$ in \eqref{cotab_univ} is determined by a (unique) morphism $E\to 1$, and by a pair of functions $A \xto p X \xot q B$ such that $\cop[s]pq : A+B\to X$ coequalizes the pair $(\inj_B\circ s,\inj_A\circ \pi_A)$.
	\end{itemize}
	(This phenomenon is similar to what happens in the double category of relations, \cite[3.6.2]{Grandis2019} where the cotabulator of a function $f : X\to Y$ is also a certain quotient of $X+Y$.)
\end{example}
Instead, a tabulator for $\fkf=\langle d,s\rangle : A \pto B$ is a cell of shape  $\footnotesize\begin{tikzcd}[cramped]
		{\tabu \fkf} & {\tabu \fkf} \\
		A & B
		\arrow["\shortmid"{marking}, Rightarrow, no head, from=1-1, to=1-2]
		\arrow["a"', from=1-1, to=2-1]
		\arrow["b", from=1-2, to=2-2]
		\arrow["{\varpi_\fkf}"{description}, draw=none, from=2-1, to=1-2]
		\arrow["\fkf", "\shortmid"{marking}, from=2-1, to=2-2]
	\end{tikzcd}$ subject to a universal property dual to the one above; the existence of such a cell implies that the carrier $E$ of $\fkf$ is nonempty (a condition that, as mild as it is, already rules out the existence of \emph{all} tabulators), and pointed by a certain element $e_0=\varpi_\fkf(\top)\in E$ in a universal way; the universal property can then only be satisfied by taking a suitable subobject of $A\times B$ as carrier of the tabulator $\idpto{\tabu \fkf}$; such a subobject $U$ must satisfy the cell condition of $\varpi_\fkf$ that reads as $\forall (a,b)\in U\subseteq A\times B : s (a , e) = b$, which is clearly absurd assuming $s$ must be a function. Moreover, $e_0$ must be such that every other cell $\xi$ as in (the dual of) \eqref{cotab_univ} factors through $e_0$, which is also absurd if $E$ has more than 2 elements.
\subsection{Companions and conjoints}
The other important genuinely double categorical universal that exists in $\cartDblMly$ is a companion for every tight morphism.
\begin{definition}[Companion of a tight morphism]\label{def_cumpa}
	Let $\dblD$ be a \PDC; the \emph{companion} of a tight morphism $f : A\to B$ is a loose morphism $f_* : A\pto B$ equipped with cells
	\[\begin{tikzcd}[cramped]
			A & B & A & A \\
			B & B & A & B
			\arrow["{f_*}", "\shortmid"{marking}, from=1-1, to=1-2]
			\arrow["f"', from=1-1, to=2-1]
			\arrow["{\epsilon_f}"{description}, draw=none, from=1-1, to=2-2]
			\arrow[Rightarrow, no head, from=1-2, to=2-2]
			\arrow["\shortmid"{marking}, Rightarrow, no head, from=1-3, to=1-4]
			\arrow[Rightarrow, no head, from=1-3, to=2-3]
			\arrow["{\eta_f}"{description}, draw=none, from=1-3, to=2-4]
			\arrow["f", from=1-4, to=2-4]
			\arrow["\shortmid"{marking}, Rightarrow, no head, from=2-1, to=2-2]
			\arrow["{f_*}"', "\shortmid"{marking}, from=2-3, to=2-4]
		\end{tikzcd}\]
	subject to the companion identities:
	\begin{itemize}
		\item the cell $\epsilon_f \vertC \eta_f$ equals the identity of $f$;
		\item the cell $\epsilon_f \horC \eta_f$ equals the identity of $f_*$.
	\end{itemize}
\end{definition}
Companions are defined by a universal property, from which they are determined uniquely (hence the fact that we call it `the' companion of $f$), that, however, we do not need to invoke.

Dually, the conjoint of a tight morphism $f : A \to B$ is the companion of $f^\vop : B\to A$ in $\dblD^\vop$.
\begin{proposition}\label{has_cumpa}
	(Proof below) Every tight morphism in $\cartDblMly$ has a companion.
\end{proposition}
\begin{example}
	A tight morphism in $\cartDblMly$ has a conjoint if and only if it is an isomorphism; give such a tight morphism $f : A \to B$ one has to build a loose morphism $\fkf^* : X \xot{d_{f^*}} B\times X \xto{s_{f^*}} A$, equipped with cells of sorts, forcing the carrier $X$ of $\fkf^*$ to be pointed by $e_0$, and the equations
	\[s (f a , e_0) = a
		\quad
		d (f a , e_0) = e_0
		\quad
		b = f (s (b , e))
	\]
	to be true; regardless of the fact that these force $f$ to be a map of a special sort (the last: $f$ is a retraction of all $s(\firstblank,e)$'s at the same time; the first: $f$ is also a retraction of $s(\firstblank,e_0)$), it's one of the cell equalities that forces the carrier of $\fkf^*$ to be a singleton, so that $\fkf^*$ is the companion of a tight morphism $f' : B\to A$; all these conditions together, however, yield that $f$ is an isomorphism with inverse $f' : B\cong B\times 1\to A$.
\end{example}
Fortunately, companions are still sufficient to link the universal Mealy fibration with the \PDC $\cartDblMly$ studied here:  we expand on the matter, leveraging on the results in \cite{lawler2015fibrationspredicatesbicategoriesrelations}, in \autoref{gnogni} below. In particular, \cite[§3.1]{loregian2024automata} observes that one can package the universal Mealy fibration into a 2-promonad, and \cite{lawler2015fibrationspredicatesbicategoriesrelations} observes that this is sufficient to build a proarrow equipment \cite{cruttwell2010unified,rosebrugh1988proarrows,wood1982abstract,wood1985proarrows} out of the (functorial) assignment $\clV(\cartDblMly) \to \clH(\cartDblMly) : f\mapsto f_*$ sending a tight morphism to its companion.
\section{From the universal fibration to the \PDC of Mealy automata, and back}\label{gnogni}
Recall from \cite{loregian2024automata} the following construction: given a category $\clK$, a diagram $B : \clB\to\clK$, and an endofunctor $F : \clK\to\clK$, the category of $F$-Mealy automata $\Mly(F,B)$ is obtained as follows:\footnote{A consequence of this description, which we will not exploit here but which will be the subject of future analysis, is that $\Mly(F,B)$ is the category of models of a $\Cat$-enriched sketch, in the sense of \cite{borceux1996enriched,borceux1998theory,Borceux2008}.}
\begin{itemize}
	\item take the \emph{inserter} $\Alg(F)$ of $F$ and $\id_\clK$, to build the object of $F$-algebras; then take the \emph{comma object} $F/B$;
	      \[\begin{tikzcd}[cramped,sep=small]
			      & {\Alg(F)} && {F/B} \\
			      \clK && \clK && \clB \\
			      & \clK && \clK
			      \arrow["U"', from=1-2, to=2-1]
			      \arrow["U", from=1-2, to=2-3]
			      \arrow["V"', from=1-4, to=2-3]
			      \arrow["{U'}", from=1-4, to=2-5]
			      \arrow[""{name=0, anchor=center, inner sep=0}, "F"', from=2-1, to=3-2]
			      \arrow[""{name=1, anchor=center, inner sep=0}, Rightarrow, no head, from=2-3, to=3-2]
			      \arrow[""{name=2, anchor=center, inner sep=0}, "F"', from=2-3, to=3-4]
			      \arrow[""{name=3, anchor=center, inner sep=0}, "B", from=2-5, to=3-4]
		      \end{tikzcd}\]
	\item now take the strict pullback of the cospan $U,V$; this is $\Mly(F,B)$, which now has been defined in a fashion that suits to any 2-category with finite limits.

	      The objects and morphisms of $\Mly(F,B)$ are, respectively, obtained as spans $\langle X,Y,d,s\rangle : X \xot d FX \xto s BY$ and pairs $u,v$ of morphisms in $\clK$ such that the squares in
	      \[\begin{tikzcd}[cramped]
			      X & FX & BY \\
			      {X'} & {FX'} & {BY'}
			      \arrow["u"', from=1-1, to=2-1]
			      \arrow["d"', from=1-2, to=1-1]
			      \arrow["s", from=1-2, to=1-3]
			      \arrow["Fu"', from=1-2, to=2-2]
			      \arrow["Bv", from=1-3, to=2-3]
			      \arrow["{d'}"', from=2-2, to=2-1]
			      \arrow["{s'}", from=2-2, to=2-3]
		      \end{tikzcd}\]
	      both commute.
\end{itemize}
In particular, when $\clK$ is Cartesian, $B : 1\to\clK$ is an object of $\clK$, and $F = A\times\firstblank$ is the functor `product by $A\in\clK$', we get back the category $\Mly(A,B)$ whose objects are spans like \eqref{florpo}. This can be conveniently packaged in the presence of a pseudofunctor
\[\begin{tikzcd}[cramped]
		{\Mly:\clK^\op\times\clK} & \Cat
		\arrow[from=1-1, to=1-2]
	\end{tikzcd}\]
in turn giving rise to a split Grothendieck fibration over $\clK^\op\times\clK$ called the \emph{universal Mealy fibration} in \cite{loregian2024automata}. In this section we describe the relation in which the universal Mealy fibration $\Mly$ and the \PDC $\dblMly$ built in \autoref{def_dbl_processes} stand.

The idea is to exploit the result of \cite{lawler2015fibrationspredicatesbicategoriesrelations}, where it is hown that the following pieces of data are equivalent:
\begin{itemize}
	\item a $\Cat$-enriched promonad $M : \clK^\op\times\clK\to\Cat$, i.e. a pseudofunctor so that, regarded as a pseudoprofunctor $\clK\pro \clK$, satisfies the axioms of a pseudomonad inside the 3-category of 2-categories, pseudoprofunctors, and pseudoprofunctor 2- and 3-cells.
	\item A double category $\dblCat{D}$ admitting all companions.
\end{itemize}
This construction makes the functor $\dblCat{D}_t \to \dblCat{D}_\ell$ sending a tight 1-cell to its companion a proarrow equipment; the loose bicategory $\dblCat{D}_\ell$ boils down to be the Kleisli bicategory of such an identity-on-objects map.

\section{Future work}
In the case of a non-Cartesian monoidal category $\clK$ of processes some aspects of the theory get more interesting; there are nonCartesian monoidal structures that are still manageable and can be described, e.g. products of convolution type, when $\clK$ is a category of (possibly enriched) presheaves (e.g. when $\clK$ is the category of $R$-modules, or a category of representations of a groupoid like in \cite{loregian2024automata}); or, semiCartesian categories given by taking as $\clK$ the Kleisli category $\Kl(D)$ of a probability monad $D$ \cite{Fritz2020,Jacobs2018,Doberkat2007,Fritz2023}. It would be interesting to study the semiCartesian \PDC of probabilistic Mealy automata.

Monads also become more interesting, and comonads nontrivial, in situations where the monoidal structure is canonical enough.

The co/completeness properties of $\cartDblMly$ make it a strange animal: tight colimit fail to exist for a trivial reason; tight limits exist and are induced by $\clK$, but \emph{double} limits fail to exist due to the absence of tabulators. The existence of some universal objects, like Kan extensions, Kan lifts, and co/comma objects might be interesting to rule out, or characterize (e.g. with a colimit formula). The search for co/completeness theorems is not only motivated by purely category-theoretic interest; often, we can give computational meaning to universal properties, and constructions that are nicely described in the bicategory probably acquire a better description in the double category: for example, behaviour as Kan extension \cite{Bainbridge1975}, or minimization in the sense of \cite{lmcs_6213}.

Tackling the former problem goes in the direction of further generalization: in \cite{Grandis2019} the authors define an \emph{internal hom} double category $\dblCat{[C,D]}$ in the category of \PDC{s} (the definition is space consuming and it is recalled below, cf. \autoref{exp_o_dbl}); replacing Walters' bicategory of pseudofunctors $\bbN\to\Sigma\clK$ with the double category of pseudo double functors from $\bbN$ (having discrete tight category, and the bicategory $\Sigma\bbN$ as loose bicategory) into $\Sigma\clK$ should lead to some insight about its universal properties \emph{qua} \PDC.

\bibliography{allofthem}{}
\bibliographystyle{plain}

\appendix
\section{Paralipomena and proofs}
\begin{definition}[\PDC]\label{def_pdc}
	A \emph{\PDC} consists of
	\begin{enumerate}
		\item \label{def_dbl_1} a class of objects $X,Y,Z,\dots$;
		\item \label{def_dbl_2} a category $\clD_0$ of \emph{tight morphisms};
		\item \label{def_dbl_3} for every pair of objects $X,Y$ a set of \emph{loose morphisms} $p : X\pto Y$;
		\item \label{def_dbl_4} for every quadruple of tight and loose morphisms sharing boundary objects as in
		      \[\begin{tikzcd}[cramped]
				      X & Y \\
				      Z & W
				      \arrow[""{name=0, anchor=center, inner sep=0}, "p", "\shortmid"{marking}, from=1-1, to=1-2]
				      \arrow["f"', from=1-1, to=2-1]
				      \arrow["g", from=1-2, to=2-2]
				      \arrow[""{name=1, anchor=center, inner sep=0}, "q"', "\shortmid"{marking}, from=2-1, to=2-2]
			      \end{tikzcd}\]
		      a set of \emph{cells} $\alpha$ having $p$ as \emph{loose domain}, $q$ as \emph{loose codomain}, $f$ as \emph{tight domain}, $g$ as \emph{tight codomain};
		\item \label{def_dbl_5} distinguished loose morphisms denoted $i_X : \idpto X$, and distinguished cells
		      \[\begin{tikzcd}[cramped]
				      A & A \\
				      B & B
				      \arrow["{i_A}",Rightarrow, no head,  "\shortmid"{marking}, from=1-1, to=1-2]
				      \arrow["f"', from=1-1, to=2-1]
				      \arrow["{\iota_f}"{description}, draw=none, from=1-1, to=2-2]
				      \arrow["f", from=1-2, to=2-2]
				      \arrow["{i_B}"', Rightarrow, no head, "\shortmid"{marking}, from=2-1, to=2-2]
			      \end{tikzcd}\]
		      for every tight morphism $f : A \to B$, and distinguished cells
		      \[\begin{tikzcd}[cramped]
				      X & Y \\
				      X & Y
				      \arrow["p", "\shortmid"{marking}, from=1-1, to=1-2]
				      \arrow[Rightarrow, no head, from=1-1, to=2-1]
				      \arrow["{\iota_p}"{description}, draw=none, from=1-1, to=2-2]
				      \arrow[Rightarrow, no head, from=1-2, to=2-2]
				      \arrow["p"', "\shortmid"{marking}, from=2-1, to=2-2]
			      \end{tikzcd}\]
		      for every loose morphism $p : X\to Y$ (where vertically we put the identities of $X,Y$ in the category $\clD_0$);
		\item \label{def_dbl_6} composition laws for cells, along the loose direction,
		      \[\begin{tikzcd}[cramped]
				      X & Y & Z && X & Z \\
				      {X'} & {Y'} & {Z'} && {X'} & {Z'}
				      \arrow["p", "\shortmid"{marking}, from=1-1, to=1-2]
				      \arrow["f"', from=1-1, to=2-1]
				      \arrow["\alpha"{description}, draw=none, from=1-1, to=2-2]
				      \arrow["{p'}", "\shortmid"{marking}, from=1-2, to=1-3]
				      \arrow["g"', from=1-2, to=2-2]
				      \arrow["\beta"{description}, draw=none, from=1-2, to=2-3]
				      \arrow["h", from=1-3, to=2-3]
				      \arrow["\mapsto"{description}, draw=none, from=1-3, to=2-5]
				      \arrow["{p'\horC p}", "\shortmid"{marking}, from=1-5, to=1-6]
				      \arrow["f"', from=1-5, to=2-5]
				      \arrow["{\beta\horC\alpha}"{description}, draw=none, from=1-5, to=2-6]
				      \arrow["h", from=1-6, to=2-6]
				      \arrow["q"', "\shortmid"{marking}, from=2-1, to=2-2]
				      \arrow["{q'}"', "\shortmid"{marking}, from=2-2, to=2-3]
				      \arrow["{q'\horC q}"', "\shortmid"{marking}, from=2-5, to=2-6]
			      \end{tikzcd}\]
		      and along the tight direction, here prescribed according to the category structure on $\clD_0$:
		      \[\begin{tikzcd}[row sep=4mm,cramped]
				      X & A && X & A \\
				      Y & B \\
				      Z & C && Z & C
				      \arrow["p", "\shortmid"{marking}, from=1-1, to=1-2]
				      \arrow["f"', from=1-1, to=2-1]
				      \arrow["\alpha"{description}, draw=none, from=1-1, to=2-2]
				      \arrow["g", from=1-2, to=2-2]
				      \arrow["\mapsto"{description}, draw=none, from=1-2, to=3-4]
				      \arrow["p", "\shortmid"{marking}, from=1-4, to=1-5]
				      \arrow["{f'f}"', from=1-4, to=3-4]
				      \arrow["{\beta\vertC\alpha}"{description}, draw=none, from=1-4, to=3-5]
				      \arrow["{g'g}", from=1-5, to=3-5]
				      \arrow["q"', "\shortmid"{marking}, from=2-1, to=2-2]
				      \arrow["{f'}"', from=2-1, to=3-1]
				      \arrow["\beta"{description}, draw=none, from=2-1, to=3-2]
				      \arrow["{g'}", from=2-2, to=3-2]
				      \arrow["r"', "\shortmid"{marking}, from=3-1, to=3-2]
				      \arrow["r"', "\shortmid"{marking}, from=3-4, to=3-5]
			      \end{tikzcd}\]
		\item \label{def_dbl_7} \emph{associator} and \emph{unitor} invertible cells
		      \[\begin{tikzcd}[cramped]
				      X && Z & W & X & Y & Y & X & X & Y \\
				      X & Y && W & X && Y & X && Y.
				      \arrow["{q\horC r}", "\shortmid"{marking}, from=1-1, to=1-3]
				      \arrow[Rightarrow, no head, from=1-1, to=2-1]
				      \arrow["{\mathrm{assoc}}"{description}, draw=none, from=1-1, to=2-4]
				      \arrow["p", "\shortmid"{marking}, from=1-3, to=1-4]
				      \arrow[Rightarrow, no head, from=1-4, to=2-4]
				      \arrow["p", "\shortmid"{marking}, from=1-5, to=1-6]
				      \arrow[Rightarrow, no head, from=1-5, to=2-5]
				      \arrow["{\mathrm{unit}^L}"{description}, draw=none, from=1-5, to=2-7]
				      \arrow["{i_Y}", "\shortmid"{marking}, from=1-6, to=1-7]
				      \arrow[Rightarrow, no head, from=1-7, to=2-7]
				      \arrow["{i_X}", "\shortmid"{marking}, from=1-8, to=1-9]
				      \arrow[Rightarrow, no head, from=1-8, to=2-8]
				      \arrow["{\mathrm{unit}^R}"{description}, draw=none, from=1-8, to=2-10]
				      \arrow["p", "\shortmid"{marking}, from=1-9, to=1-10]
				      \arrow[Rightarrow, no head, from=1-10, to=2-10]
				      \arrow["r"', "\shortmid"{marking}, from=2-1, to=2-2]
				      \arrow["{p\horC q}"', "\shortmid"{marking}, from=2-2, to=2-4]
				      \arrow["p"', "\shortmid"{marking}, from=2-5, to=2-7]
				      \arrow["p"', "\shortmid"{marking}, from=2-8, to=2-10]
			      \end{tikzcd}\]
	\end{enumerate}
	All these data are subject to the conditions stated in \cite[§7.1]{CTGDC_1999__40_3_162_0}, particularly, [\emph{ibi}, pdc2.1---pdc2.6]:
	\begin{enumerate}
		\item \label{da_1} (weak) associativity of loose composition: given composable cells $\alpha,\beta,\gamma$ the two possible loose compositions $\alpha\mid \gamma\horC\beta$ and $\beta\horC\alpha\mid\gamma$;
		      differ by the associator in \autoref{def_dbl_7};
		\item \label{da_2} (weak) unitality of loose composition;
		\item \label{da_2'} (weak) associativity of loose composition;
		\item bicategory axioms (pentagon and triangle) for loose composition;
		\item \label{da_3} axioms of category for the tight composition;
		\item \label{da_4} interchange law.
	\end{enumerate}
\end{definition}
The adjective `weak' is what qualifies the resulting structure as a \emph{pseudo} double category. One situation when a double category arises naturally is when a category $\clC$ is given, together with a bicategory having the same objects and a natural choice of cells as in \eqref{baiadiai}. This is for example the case of the double category of modules, cf. \cite[Example 2.2]{Shulman2008a}, the more general double category of distributors, and the double category of relations in a regular category, cf. \cite{dawson2010span} and \cite{Freyd1990}.
\begin{definition}[Duals of a \PDC]\label{los_duales}
	There are three notions of dual for a \PDC $\dblD$:
	\begin{itemize}
		\item the loose dual $\dblD^\hop$ defined reversing the direction of the loose morphisms;
		\item the tight dual $\dblD^\vop$ defined reversing the direction of the tight morphisms;
		\item the double dual $\dblD^\oop$ defined reversing the direction of both loose and tight morphisms.
	\end{itemize}
\end{definition}
\begin{definition}
	The cotabulator of a loose morphism $F : A\pto B$ in a \PDC $\dblD$ consists of a pair $(\cotabu F, \tau_F)$ where $\cotabu F$ is an object of $\dblD$, and $\tau_F$ a cell
	$\begin{tikzcd}[cramped]
			A & B \\
			{\cotabu F} & {\cotabu F}
			\arrow["F", "\shortmid"{marking}, from=1-1, to=1-2]
			\arrow["a"', from=1-1, to=2-1]
			\arrow["{\tau_F}"{description}, draw=none, from=1-1, to=2-2]
			\arrow["b", from=1-2, to=2-2]
			\arrow["\shortmid"{marking}, Rightarrow, no head, from=2-1, to=2-2]
		\end{tikzcd}$
	universal among all such, i.e. with the property that any other cell
	$\begin{tikzcd}[cramped]
			A & B \\
			X & X
			\arrow["F", "\shortmid"{marking}, from=1-1, to=1-2]
			\arrow["p"', from=1-1, to=2-1]
			\arrow["\xi"{description}, draw=none, from=1-1, to=2-2]
			\arrow["q", from=1-2, to=2-2]
			\arrow["\shortmid"{marking}, Rightarrow, no head, from=2-1, to=2-2]
		\end{tikzcd}$
	factors uniquely through $\tau_F$ as
	\[\label{cotab_univ}\begin{tikzcd}[row sep=4mm,cramped]
			A & B & A & B \\
			&& {\cotabu F} & {\cotabu F} \\
			X & X & X & X
			\arrow["F", "\shortmid"{marking}, from=1-1, to=1-2]
			\arrow["p"', from=1-1, to=3-1]
			\arrow["\xi"{description}, draw=none, from=1-1, to=3-2]
			\arrow["q", from=1-2, to=3-2]
			\arrow["{=}"{description}, draw=none, from=1-2, to=3-3]
			\arrow["F", "\shortmid"{marking}, from=1-3, to=1-4]
			\arrow[from=1-3, to=2-3]
			\arrow["{\tau_F}"{description}, draw=none, from=1-3, to=2-4]
			\arrow[from=1-4, to=2-4]
			\arrow["\shortmid"{marking}, Rightarrow, no head, from=2-3, to=2-4]
			\arrow[from=2-3, to=3-3]
			\arrow["{\hat\xi}"{description}, draw=none, from=2-3, to=3-4]
			\arrow[from=2-4, to=3-4]
			\arrow["\shortmid"{marking}, Rightarrow, no head, from=3-1, to=3-2]
			\arrow["\shortmid"{marking}, Rightarrow, no head, from=3-3, to=3-4]
		\end{tikzcd}\]
\end{definition}
The \emph{tabulator} of a loose morphism $F : A\pto B$ is defined as the cotabulator of $F$ in $\dblD^\vop$, and it consists of a universal cell
$\begin{tikzcd}[cramped]
		{\tabu F} & {\tabu F} \\
		A & B
		\arrow["\shortmid"{marking}, Rightarrow, no head, from=1-1, to=1-2]
		\arrow["a"', from=1-1, to=2-1]
		\arrow["b", from=1-2, to=2-2]
		\arrow["{\varpi_F}"{description}, draw=none, from=2-1, to=1-2]
		\arrow["F", "\shortmid"{marking}, from=2-1, to=2-2]
	\end{tikzcd}$
subject to the dual universal property of\eqref{cotab_univ}.

In the following definitions we outline the notion of (tight) co/limit in a \PDC $\dblD$, drawing from \cite[§4.2]{CTGDC_1999__40_3_162_0}. We do not need to describe colimits explicitly, as they seldom exist in $\cartDblMly$, so we are content with the tautological definition of `limits in the vertical dual of $\dblD$'.
\begin{definition}[Initial and terminal objects in a \PDC]\label{def_initerm_pc}
	A \emph{terminal object} in $\dblD$ consists of an object $\top$ which is terminal in the tight category $\clD_0$, and such that for every loose morphism $p : A\pto B$ there is a unique cell
	\[\begin{tikzcd}
			A \ar[d, "!_A"']\ar[r, "\shortmid"{marking}, "p"] & B\ar[d, "!_B"]\\
			\top \ar[r, "\shortmid"{marking}, Rightarrow,no head]& \top .
			\arrow["{!_p}"{description}, draw=none, from=1-1, to=2-2]
		\end{tikzcd}\]
	Dually, an initial object in $\dblD$ consists of a terminal object in $\dblD^\vop$.
\end{definition}
(Generalizing the following definition to the case of a family $(A_i\mid i\in I)$ of objects is straightforward.)
\begin{definition}[Products and coproducts in a \PDC]\label{def_procopro_pc}
	Let $A,B$ be two objects of $\dblD$; the \emph{double product} $A\times B$ of the two objects consists of a pair of cells
	\[\begin{tikzcd}[cramped]
			{A\times B} & {A\times B} & {A\times B} & {A\times B} \\
			A & A & B & B
			\arrow["\shortmid"{marking}, Rightarrow, no head, from=1-1, to=1-2]
			\arrow["{p_A}"', from=1-1, to=2-1]
			\arrow["{\pi_A}"{description}, draw=none, from=1-1, to=2-2]
			\arrow["{p_A}", from=1-2, to=2-2]
			\arrow["\shortmid"{marking}, Rightarrow, no head, from=1-3, to=1-4]
			\arrow["{p_B}"', from=1-3, to=2-3]
			\arrow["{\pi_B}"{description}, draw=none, from=1-3, to=2-4]
			\arrow["{p_B}", from=1-4, to=2-4]
			\arrow["\shortmid"{marking}, Rightarrow, no head, from=2-1, to=2-2]
			\arrow["\shortmid"{marking}, Rightarrow, no head, from=2-3, to=2-4]
		\end{tikzcd}\]
	that are universal in the following sense: given any loose morphism $u : X\pto Y$, tight morphisms $A \xot a X \xto b B$ and $A \xot {a'} Y \xto {b'} B$ and cells
	\[\begin{tikzcd}[cramped]
			X & Y & X & Y \\
			A & A & B & B
			\arrow["u", "\shortmid"{marking}, from=1-1, to=1-2]
			\arrow["a"', from=1-1, to=2-1]
			\arrow["{\xi_1}"{description}, draw=none, from=1-1, to=2-2]
			\arrow["{a'}", from=1-2, to=2-2]
			\arrow["u", "\shortmid"{marking}, from=1-3, to=1-4]
			\arrow["b"', from=1-3, to=2-3]
			\arrow["{\xi_2}"{description}, draw=none, from=1-3, to=2-4]
			\arrow["{b'}", from=1-4, to=2-4]
			\arrow["\shortmid"{marking}, Rightarrow, no head, from=2-1, to=2-2]
			\arrow["\shortmid"{marking}, Rightarrow, no head, from=2-3, to=2-4]
		\end{tikzcd}\]
	\begin{enumerate}
		\item there exists a unique pair of tight morphisms $\langle a,b\rangle : X\to A\times B$ and $\langle a',b'\rangle : Y\to A\times B$ so that the tuple $(A\times B; p_A,p_B)$ is a product in $\clD_0$ (so, in particular, $p_A\circ\langle a,b\rangle = a$, $p_B\circ \langle a,b\rangle = b$, $p_A\circ \langle a',b'\rangle = a'$, $p_B\circ \langle a',b'\rangle = b'$);
		\item there exists a unique cell
		      \[\begin{tikzcd}[cramped]
				      X & Y \\
				      {A\times B} & {A\times B}
				      \arrow["u", "\shortmid"{marking}, from=1-1, to=1-2]
				      \arrow["{\langle a,b\rangle}"', from=1-1, to=2-1]
				      \arrow["{\langle\xi_1,\xi_2\rangle}"{description}, draw=none, from=1-1, to=2-2]
				      \arrow["{\langle a',b'\rangle}", from=1-2, to=2-2]
				      \arrow["\shortmid"{marking}, Rightarrow, no head, from=2-1, to=2-2]
			      \end{tikzcd}\]
		      such that the equalities of cells $\xi_1 = \frac{\langle \xi_1,\xi_2\rangle}{\pi_A}$ and $\xi_2 = \frac{\langle \xi_1,\xi_2\rangle}{\pi_B}$ hold.
	\end{enumerate}
	The \emph{double coproduct} of $A, B$ consists of the double product of $A,B$ in $\dblD^\vop$.
\end{definition}
\begin{definition}[Pullbacks and pushouts in a \PDC]\label{def_pullpush_pc}
	Let $A \xto f X \xot g B$ be a cospan of tight morphisms in $\dblD$; the \emph{double pullback} of $f,g$ in $\dblD$ consists of an arrangement of cells
	\[\begin{tikzcd}[cramped]
			{A\times_X B} & {A\times_X B} & {A\times_X B} & {A\times_X B} \\
			A & A & B & B
			\arrow["\shortmid"{marking}, Rightarrow, no head, from=1-1, to=1-2]
			\arrow["{p_A}"', from=1-1, to=2-1]
			\arrow["{\pi_A}"{description}, draw=none, from=1-1, to=2-2]
			\arrow["{p_A}", from=1-2, to=2-2]
			\arrow["\shortmid"{marking}, Rightarrow, no head, from=1-3, to=1-4]
			\arrow["{p_B}"', from=1-3, to=2-3]
			\arrow["{\pi_B}"{description}, draw=none, from=1-3, to=2-4]
			\arrow["{p_B}", from=1-4, to=2-4]
			\arrow["\shortmid"{marking}, Rightarrow, no head, from=2-1, to=2-2]
			\arrow["\shortmid"{marking}, Rightarrow, no head, from=2-3, to=2-4]
		\end{tikzcd}\]
	so that $f\circ p_A = g\circ p_B$, and terminal among all such. This means that
	given any loose morphism $u : U\pto V$, tight morphisms $A \xot a U \xto b B$ and $A \xot {a'} V \xto {b'} B$ and cells
	\[\begin{tikzcd}[cramped]
			U & V & U & V \\
			A & A & B & B
			\arrow["u", "\shortmid"{marking}, from=1-1, to=1-2]
			\arrow["a"', from=1-1, to=2-1]
			\arrow["{\xi_A}"{description}, draw=none, from=1-1, to=2-2]
			\arrow["{a'}", from=1-2, to=2-2]
			\arrow["u", "\shortmid"{marking}, from=1-3, to=1-4]
			\arrow["b"', from=1-3, to=2-3]
			\arrow["{\xi_B}"{description}, draw=none, from=1-3, to=2-4]
			\arrow["{b'}", from=1-4, to=2-4]
			\arrow["\shortmid"{marking}, Rightarrow, no head, from=2-1, to=2-2]
			\arrow["\shortmid"{marking}, Rightarrow, no head, from=2-3, to=2-4]
		\end{tikzcd}\]
	such that $\frac{\xi_A}f = \frac{\xi_B}g$ (so in particular, $f\circ a = g\circ b$ and $f\circ a' = g\circ b'$),
	\begin{enumerate}
		\item there exists a unique pair of tight morphisms $\langle a,b\rangle : U\to A\times_X B$ and $\langle a',b'\rangle : V\to A\times_X B$, so that in particular $A\times_X B$ is a pullback in $\clD_0$;
		\item there exists a unique factorization of cells $\xi_A,\xi_B$ as 
		      \[\begin{tikzcd}[cramped]
				      & U \\
				      & {A\times_X B} & A & V \\
				      B &&& {A\times_X B} & A \\
				      & X & B \\
				      &&& X
				      \arrow["{!}", dashed, from=1-2, to=2-2]
				      \arrow[from=1-2, to=2-3]
				      \arrow["u", "\shortmid"{marking}, from=1-2, to=2-4]
				      \arrow[curve={height=12pt}, from=1-2, to=3-1]
				      \arrow[from=2-2, to=2-3]
				      \arrow["{p_B}"{description}, from=2-2, to=3-1]
				      \arrow["\shortmid"{marking}, Rightarrow, no head, from=2-2, to=3-4]
				      \arrow["\lrcorner"{anchor=center, pos=0.125, rotate=-45}, draw=none, from=2-2, to=4-2]
				      \arrow["\shortmid"{marking}, Rightarrow, no head, from=2-3, to=3-5]
				      \arrow["f"', from=2-3, to=4-2]
				      \arrow["{!}"'{pos=0.8}, dashed, from=2-4, to=3-4]
				      \arrow[from=2-4, to=3-5]
				      \arrow[curve={height=12pt}, from=2-4, to=4-3]
				      \arrow["g"', from=3-1, to=4-2]
				      \arrow["\shortmid"{marking}, Rightarrow, no head, from=3-1, to=4-3]
				      \arrow["{p_A}"', from=3-4, to=3-5]
				      \arrow["{p_B}"{description}, from=3-4, to=4-3]
				      \arrow["\lrcorner"{anchor=center, pos=0.125, rotate=-45}, draw=none, from=3-4, to=5-4]
				      \arrow["f", from=3-5, to=5-4]
				      \arrow["\shortmid"{marking}, Rightarrow, no head, from=4-2, to=5-4]
				      \arrow["g", from=4-3, to=5-4]
			      \end{tikzcd}\]
		      through a unique mediating cell $\smallCell{u}{!}{!}{}[\langle \xi_A,\xi_B\rangle]$ such that
		      such that the equalities of cells $\xi_A = \frac{\langle \xi_A,\xi_B\rangle}{\pi_A}$ and $\xi_B = \frac{\langle \xi_A,\xi_B\rangle}{\pi_B}$ hold.
	\end{enumerate}
	The \emph{double pushout} of $A, B$ consists of the double product of $A,B$ in $\dblD^\vop$.
\end{definition}
\begin{definition}[Exponential of double categories]\label{exp_o_dbl}
	Let $\dblC,\dblD$ be two double categories; we construct the internal hom double category $\dblCat{[C,D]}$ as follows.
	\begin{itemize}
		\item The objects are pseudo double functors\footnote{A version $\dblCat{[C,D]}_\ell$ of this construction exists where objects are \emph{lax} double functors and the rest of the structure is unchanged --just, tight, loose morphisms and cells satisfy more coherence axioms. We do not need such generality as Walters et al. consider only \emph{pseudo} functors $\bbN \to\Sigma\clK$ in \cite{Katis1997}.} $F : \dblC \to\dblD$, as defined in \cite[3.5.1 and p. 145]{Grandis2019};
		\item tight morphisms $\phi : F \To G$ consist of families of tight morphisms $\phi_C : FC\to GC$ in $\dblD$ which are \emph{tightly natural}, i.e. for every tight $u : C\to C'$ one has $\phi_{C'}\circ Fu = Gu\circ \phi_C$, and \emph{loosely lax natural}, i.e. for every loose $c : C\pto C'$ there are cells $\smallCell{Fc}{\phi_C}{\phi_{C'}}{Gc}[\phi_c]$ such that
		      \begin{itemize}
			      \item the $\phi_c$ are compatible with the structure isomorphisms of $F,G$, i.e. with the loose isomorphisms $Fu\circ Fv \cong F(u\circ v)$, $\iota_{FC}\cong F(\iota_C)$, $Gu\circ Gv \cong G(u\circ v)$, $\iota_{GC}\cong G(\iota_C)$: for example, any time $c'\odot_h c : C \pto C' \pto C''$ are composable,
			            \[\begin{tikzcd}[cramped]
					            FC & {FC'} & {FC''} && FC & {FC'} & {FC''} \\
					            GC & {GC'} & {GC''} && FC && {FC''} \\
					            GC && {GC''} && GC && {GC''}
					            \arrow["Fc", "\shortmid"{marking}, from=1-1, to=1-2]
					            \arrow["{\phi_C}"', from=1-1, to=2-1]
					            \arrow["{\phi_c}"{description}, draw=none, from=1-1, to=2-2]
					            \arrow["{Fc'}", "\shortmid"{marking}, from=1-2, to=1-3]
					            \arrow[from=1-2, to=2-2]
					            \arrow["{\phi_{c'}}"{description}, draw=none, from=1-2, to=2-3]
					            \arrow["{\phi_{C''}}", from=1-3, to=2-3]
					            \arrow["{=}"{description}, draw=none, from=1-3, to=3-5]
					            \arrow["Fc", "\shortmid"{marking}, from=1-5, to=1-6]
					            \arrow[Rightarrow, no head, from=1-5, to=2-5]
					            \arrow["\cong"{description}, draw=none, from=1-5, to=2-7]
					            \arrow["{Fc'}", "\shortmid"{marking}, from=1-6, to=1-7]
					            \arrow[Rightarrow, no head, from=1-7, to=2-7]
					            \arrow["\shortmid"{marking}, from=2-1, to=2-2]
					            \arrow[Rightarrow, no head, from=2-1, to=3-1]
					            \arrow["\cong"{description}, draw=none, from=2-1, to=3-3]
					            \arrow["\shortmid"{marking}, from=2-2, to=2-3]
					            \arrow[Rightarrow, no head, from=2-3, to=3-3]
					            \arrow["\shortmid"{marking}, from=2-5, to=2-7]
					            \arrow["{\phi_C}"', from=2-5, to=3-5]
					            \arrow["{\phi_{c'\odot_h c}}"{description}, draw=none, from=2-5, to=3-7]
					            \arrow["{\phi_{C''}}", from=2-7, to=3-7]
					            \arrow["{G(c'\odot_h c)}"', "\shortmid"{marking}, from=3-1, to=3-3]
					            \arrow["{G(c'\odot_h c)}"', "\shortmid"{marking}, from=3-5, to=3-7]
				            \end{tikzcd}\]
			      \item the $\phi_c$ are compatible with cells of $\dblD$, i.e. for every cell $\smallCell suvt[\zeta]$,
			            \[\begin{tikzcd}[cramped]
					            {\cdot} & {\cdot} && {\cdot} & {\cdot} \\
					            {\cdot} & {\cdot} & \equiv & {\cdot} & {\cdot} \\
					            {\cdot} & {\cdot} && {\cdot} & {\cdot}
					            \arrow[""{name=0, anchor=center, inner sep=0}, "Fs", "\shortmid"{marking}, from=1-1, to=1-2]
					            \arrow["Fu"', from=1-1, to=2-1]
					            \arrow["Fv", from=1-2, to=2-2]
					            \arrow[""{name=1, anchor=center, inner sep=0}, "Fs", "\shortmid"{marking}, from=1-4, to=1-5]
					            \arrow["{\phi_{\#}}"', from=1-4, to=2-4]
					            \arrow["{\phi_{\#}}", from=1-5, to=2-5]
					            \arrow[""{name=2, anchor=center, inner sep=0}, "Ft"', "\shortmid"{marking}, from=2-1, to=2-2]
					            \arrow["{\phi_{\#}}"', from=2-1, to=3-1]
					            \arrow["{\phi_{\#}}", from=2-2, to=3-2]
					            \arrow[""{name=3, anchor=center, inner sep=0}, "\shortmid"{marking}, from=2-4, to=2-5]
					            \arrow["Gu"', from=2-4, to=3-4]
					            \arrow["Gv", from=2-5, to=3-5]
					            \arrow[""{name=4, anchor=center, inner sep=0}, "Gt"', "\shortmid"{marking}, from=3-1, to=3-2]
					            \arrow[""{name=5, anchor=center, inner sep=0}, "Gt"', "\shortmid"{marking}, from=3-4, to=3-5]
					            \arrow["{F\zeta}"{description}, draw=none, from=0, to=2]
					            \arrow["{\phi_s}"{description}, draw=none, from=1, to=3]
					            \arrow["{\phi_t}"{description}, draw=none, from=4, to=2]
					            \arrow["{G\zeta}"{description}, draw=none, from=5, to=3]
				            \end{tikzcd}\]
		      \end{itemize}
		\item loose morphisms $\psi : F\Loose G$ consist of
		      \begin{itemize}
			      \item families of loose morphisms $\psi_C : FC\pto GC$,
			      \item    cells of shape
			            \[\begin{tikzcd}[cramped]
					            FC & GC & {GC'} \\
					            FC & {FC'} & {GC'}
					            \arrow["{\psi_C}", "\shortmid"{marking}, from=1-1, to=1-2]
					            \arrow[Rightarrow, no head, from=1-1, to=2-1]
					            \arrow["Gc", "\shortmid"{marking}, from=1-2, to=1-3]
					            \arrow["{\psi_c}"{description}, draw=none, from=1-2, to=2-2]
					            \arrow[Rightarrow, no head, from=1-3, to=2-3]
					            \arrow["Fc"', "\shortmid"{marking}, from=2-1, to=2-2]
					            \arrow["{\psi_{C'}}"', "\shortmid"{marking}, from=2-2, to=2-3]
				            \end{tikzcd}\]
			      \item for each $u : C\to C'$ in $\dblC$, cells in $\dblD$ of shape
			            \[\begin{tikzcd}[cramped]
					            FC & GC \\
					            {FC'} & {GC'}
					            \arrow[""{name=0, anchor=center, inner sep=0}, "{\psi_C}", "\shortmid"{marking}, from=1-1, to=1-2]
					            \arrow["Fu"', from=1-1, to=2-1]
					            \arrow["Gu", from=1-2, to=2-2]
					            \arrow[""{name=1, anchor=center, inner sep=0}, "{\psi_{C'}}"', "\shortmid"{marking}, from=2-1, to=2-2]
					            \arrow["{\psi_f}"{description}, draw=none, from=0, to=1]
				            \end{tikzcd}\]
		      \end{itemize}
		      such that for every cell $\smallCell sfgt[\zeta]$,
		      \[\begin{tikzcd}[cramped]
				      \cdot & \cdot & \cdot && \cdot & \cdot & \cdot \\
				      \cdot & \cdot & \cdot & \equiv & \cdot & \cdot & \cdot \\
				      \cdot & \cdot & \cdot && \cdot & \cdot & \cdot
				      \arrow[""{name=0, anchor=center, inner sep=0}, "{\psi_{\#}}", "\shortmid"{marking}, from=1-1, to=1-2]
				      \arrow["Fu"', from=1-1, to=2-1]
				      \arrow[""{name=1, anchor=center, inner sep=0}, "Gs", "\shortmid"{marking}, from=1-2, to=1-3]
				      \arrow[from=1-2, to=2-2]
				      \arrow["Gv", from=1-3, to=2-3]
				      \arrow["{\psi_{\#}}", "\shortmid"{marking}, from=1-5, to=1-6]
				      \arrow[Rightarrow, no head, from=1-5, to=2-5]
				      \arrow["Gs", "\shortmid"{marking}, from=1-6, to=1-7]
				      \arrow["{\psi_s}"{description}, draw=none, from=1-6, to=2-6]
				      \arrow[Rightarrow, no head, from=1-7, to=2-7]
				      \arrow[""{name=2, anchor=center, inner sep=0}, "{\psi_{\#}}"', "\shortmid"{marking}, from=2-1, to=2-2]
				      \arrow[Rightarrow, no head, from=2-1, to=3-1]
				      \arrow[""{name=3, anchor=center, inner sep=0}, "Gt"', "\shortmid"{marking}, from=2-2, to=2-3]
				      \arrow["{\psi_t}"{description}, draw=none, from=2-2, to=3-2]
				      \arrow[Rightarrow, no head, from=2-3, to=3-3]
				      \arrow[""{name=4, anchor=center, inner sep=0}, "Fs", "\shortmid"{marking}, from=2-5, to=2-6]
				      \arrow["Fu"', from=2-5, to=3-5]
				      \arrow[""{name=5, anchor=center, inner sep=0}, "\shortmid"{marking}, from=2-6, to=2-7]
				      \arrow[from=2-6, to=3-6]
				      \arrow["Gv", from=2-7, to=3-7]
				      \arrow["Ft"', "\shortmid"{marking}, from=3-1, to=3-2]
				      \arrow["{\psi_{\#}}"', "\shortmid"{marking}, from=3-2, to=3-3]
				      \arrow[""{name=6, anchor=center, inner sep=0}, "Ft"', "\shortmid"{marking}, from=3-5, to=3-6]
				      \arrow[""{name=7, anchor=center, inner sep=0}, "{\psi_{\#}}"', "\shortmid"{marking}, from=3-6, to=3-7]
				      \arrow["{\psi_f}"{description}, draw=none, from=0, to=2]
				      \arrow["{G\zeta}"{description}, draw=none, from=1, to=3]
				      \arrow["{F\zeta}"{description}, draw=none, from=4, to=6]
				      \arrow["{\psi_v}"{description}, draw=none, from=5, to=7]
			      \end{tikzcd}\]
		\item cells $\smallCell\psi\phi{\phi'}{\psi'}[\Delta]$ in $[\dblC,\dblD]$ consist of families $\Delta_C$ of cells of $\dblD$ indexed by the objects of $\dblC$,
		      \[\begin{tikzcd}[cramped]
				      F & G \\
				      {F'} & {G'}
				      \arrow[""{name=0, anchor=center, inner sep=0}, "\psi", "\shortmid"{marking}, from=1-1, to=1-2]
				      \arrow["\phi"', from=1-1, to=2-1]
				      \arrow["{\phi'}", from=1-2, to=2-2]
				      \arrow[""{name=1, anchor=center, inner sep=0}, "{\psi'}"', "\shortmid"{marking}, from=2-1, to=2-2]
				      \arrow["\Delta"{description}, draw=none, from=0, to=1]
			      \end{tikzcd}\]
		      such that for every loose morphism $c : C\pto C'$ and tight morphism $f : C\to C'$, the following conditions are satisfied:
		      \begin{itemize}
			      \item compatibility with the $\psi_c$:
			            \[\begin{tikzcd}[cramped]
					            FC & GC & {GC'} && FC & GC & {GC'} \\
					            {F'C} & {G'C} & {G'C'} & \equiv & FC & {FC'} & {GC'} \\
					            {F'C} & {F'C'} & {G'C'} && {F'C} & {F'C'} & {G'C'}
					            \arrow["{\psi_C}", "\shortmid"{marking}, from=1-1, to=1-2]
					            \arrow["{\phi_C}"', from=1-1, to=2-1]
					            \arrow["{\Delta_C}"{description}, draw=none, from=1-1, to=2-2]
					            \arrow["Gc", "\shortmid"{marking}, from=1-2, to=1-3]
					            \arrow["{\phi'_C}"{description}, from=1-2, to=2-2]
					            \arrow["{\phi'_c}"{description}, draw=none, from=1-2, to=2-3]
					            \arrow["{\phi'_{C'}}", from=1-3, to=2-3]
					            \arrow["{\psi_C}", "\shortmid"{marking}, from=1-5, to=1-6]
					            \arrow[Rightarrow, no head, from=1-5, to=2-5]
					            \arrow["{\psi_c}"{description}, draw=none, from=1-5, to=2-7]
					            \arrow["Gc", "\shortmid"{marking}, from=1-6, to=1-7]
					            \arrow[Rightarrow, no head, from=1-7, to=2-7]
					            \arrow["\shortmid"{marking}, from=2-1, to=2-2]
					            \arrow[Rightarrow, no head, from=2-1, to=3-1]
					            \arrow["{\psi'_c}"{description}, draw=none, from=2-1, to=3-3]
					            \arrow["\shortmid"{marking}, from=2-2, to=2-3]
					            \arrow[Rightarrow, no head, from=2-3, to=3-3]
					            \arrow["\shortmid"{marking}, from=2-5, to=2-6]
					            \arrow["{\phi_C}"', from=2-5, to=3-5]
					            \arrow["{\phi_c}"{description}, draw=none, from=2-5, to=3-6]
					            \arrow["\shortmid"{marking}, from=2-6, to=2-7]
					            \arrow["{\phi_{C'}}"{description}, from=2-6, to=3-6]
					            \arrow["{\Delta_{C'}}"{description}, draw=none, from=2-6, to=3-7]
					            \arrow["{\phi'_{C'}}", from=2-7, to=3-7]
					            \arrow["{F'c}"', "\shortmid"{marking}, from=3-1, to=3-2]
					            \arrow["{\psi'_{C'}}"', "\shortmid"{marking}, from=3-2, to=3-3]
					            \arrow["{F'c}"', "\shortmid"{marking}, from=3-5, to=3-6]
					            \arrow["{\psi'_{C'}}"', "\shortmid"{marking}, from=3-6, to=3-7]
				            \end{tikzcd}\]
			      \item compatibility with the $\psi_f$:
			            \[\begin{tikzcd}[cramped]
					            FC & GC && FC & GC \\
					            {FC'} & {GC'} & \equiv & {F'C} & {G'C} \\
					            {F'C'} & {G'C'} && {F'C'} & {G'C'}
					            \arrow["{\psi_C}", "\shortmid"{marking}, from=1-1, to=1-2]
					            \arrow["Fu"', from=1-1, to=2-1]
					            \arrow["{\psi_f}"{description}, draw=none, from=1-1, to=2-2]
					            \arrow["Gu", from=1-2, to=2-2]
					            \arrow["{\psi_C}", "\shortmid"{marking}, from=1-4, to=1-5]
					            \arrow["{\phi_C}"', from=1-4, to=2-4]
					            \arrow["{\Delta_C}"{description}, draw=none, from=1-4, to=2-5]
					            \arrow["{\phi'_C}", from=1-5, to=2-5]
					            \arrow["\shortmid"{marking}, from=2-1, to=2-2]
					            \arrow["{\phi_{C'}}"', from=2-1, to=3-1]
					            \arrow["{\Delta_{C'}}"{description}, draw=none, from=2-1, to=3-2]
					            \arrow["{\phi'_{C'}}", from=2-2, to=3-2]
					            \arrow["\shortmid"{marking}, from=2-4, to=2-5]
					            \arrow["{F'f}"', from=2-4, to=3-4]
					            \arrow["{\psi'_f}"{description}, draw=none, from=2-4, to=3-5]
					            \arrow["{G'f}", from=2-5, to=3-5]
					            \arrow["{\psi'_{C'}}"', "\shortmid"{marking}, from=3-1, to=3-2]
					            \arrow["{\psi'_{C'}}"', "\shortmid"{marking}, from=3-4, to=3-5]
				            \end{tikzcd}\]
		      \end{itemize}
	\end{itemize}
\end{definition}
\begin{definition}[Tight monad morphism]
	\label{mormonad_tig}
	Let $(M : A\pto A,\mu,\eta)$, $(N : B\pto B,\nu,\gamma)$ be two monads on objects $A,B$, in a \PDC $\dblD$; a \emph{(tight) monad morphism} consists of a pair $(f,\alpha)$ where $f : A\to B$ is a tight morphism, and $\alpha$ a cell
	\[\begin{tikzcd}[cramped]
			A & A \\
			B & B
			\arrow["M", "\shortmid"{marking}, from=1-1, to=1-2]
			\arrow["f"', from=1-1, to=2-1]
			\arrow["\alpha"{description}, draw=none, from=1-1, to=2-2]
			\arrow["f", from=1-2, to=2-2]
			\arrow["N"', "\shortmid"{marking}, from=2-1, to=2-2]
		\end{tikzcd}\]
	compatible with the monad structure of $M,N$ in the sense that $\frac\mu\alpha = \frac{\alpha\mid\alpha}{\nu}$ and $\frac\eta\alpha = \frac f\gamma$.
\end{definition}
This definition is used twice in the above discussion:
\begin{itemize}
	\item In \autoref{dgaddaaa}, every monad morphism in this sense	induces a monoid homomorphism $f\times\alpha : E\bowtie_M A^* \to B^*\bowtie_N E'$ through $f : A\to B$ and $\alpha : E\to E'$ (the carriers of $M$ and $N$ respectively).
	\item In \autoref{coro_fremonad}, a cell of type $\smallCell F{}{}N[\gamma]$ extends to a morphism of monads $\smallCell {M(F)}{}{}N[\gamma^*]$ promoting the construction outlined in \autoref{coro_fremonad} to a double functor.
\end{itemize}
\begin{definition}[Loose adjunction]\label{def_loose_adja}
	Two loose morphisms $l : A\pto B$, $r : B\pto A$ form a \emph{loose adjunction} in a \PDC $\dblD$ if they are equipped with
	\begin{itemize}
		\item two cells (called \emph{unit} and \emph{counit}),
		      \[\begin{tikzcd}[cramped]
				      A & A & B & B \\
				      A & A & B & B
				      \arrow["\shortmid"{marking}, Rightarrow, no head, from=1-1, to=1-2]
				      \arrow[Rightarrow, no head, from=1-1, to=2-1]
				      \arrow["\eta"{description}, draw=none, from=1-1, to=2-2]
				      \arrow[Rightarrow, no head, from=1-2, to=2-2]
				      \arrow["{l\horC r}", "\shortmid"{marking}, from=1-3, to=1-4]
				      \arrow[Rightarrow, no head, from=1-3, to=2-3]
				      \arrow["\epsilon"{description}, draw=none, from=1-3, to=2-4]
				      \arrow[Rightarrow, no head, from=1-4, to=2-4]
				      \arrow["{r\horC l}"', "\shortmid"{marking}, from=2-1, to=2-2]
				      \arrow["\shortmid"{marking}, Rightarrow, no head, from=2-3, to=2-4]
			      \end{tikzcd}\]
		\item satisfying the adjunction identities $l = \frac{\eta\mid l}{l\mid\epsilon}$ and $r = \frac{r\mid \eta}{\epsilon\mid r}$.
	\end{itemize}
\end{definition}
\subsection{Proofs}
\begin{savvyProof}[Proof of \autoref{crucial_thm}.]
	Unpacking the definition of algebra for the monad $M=((E,\langle d,s\rangle) : A \pto A,\eta,\mu)$ as in \autoref{sgrulbio}, when instantiated in $\cartDblMly$, we get
	\begin{enumtag}{ad}
		\item a Mealy automaton $(P,\langle \delta,\sigma\rangle) : A\times P \to P\times X$;
		\item an algebra map $\xi : E \times P \to P$ defining a cell
		\[\begin{tikzcd}[cramped, ampersand replacement=\&]
				A \& A \& X \\
				A \&\& X
				\arrow["M", "\shortmid"{marking}, from=1-1, to=1-2]
				\arrow[Rightarrow, no head, from=1-1, to=2-1]
				\arrow["\xi"{description}, draw=none, from=1-1, to=2-3]
				\arrow["P", "\shortmid"{marking}, from=1-2, to=1-3]
				\arrow[Rightarrow, no head, from=1-3, to=2-3]
				\arrow["P"', "\shortmid"{marking}, from=2-1, to=2-3]
			\end{tikzcd}\]
	\end{enumtag}
	subject to the following conditions (the notation is as in \autoref{deffa_monada}):
	\begin{enumtag}{ax}
		\item \label{ax_1} cell axiom for $\xi$, I: for all $x\in X, t\in P, e\in E$ we have $\sigma (a , \xi (e , x)) = \sigma (a \sact e , x)$;
		\item \label{ax_2} cell axiom for $\xi$, II: for all $x\in X, t\in P, e\in E$ we have $\delta (a , \xi (e , x)) = \xi (a \dact e , \delta (a \sact e , x))$;
		\item \label{ax_3} unitality: for all $t\in P$ we have $\xi (e_0 , x) = x$ 
		\item \label{ax_4} associativity: for $t\in P$, $e,e'\in E$ we have $\xi (\mu (e , e') , x) = \xi (e , \xi (e' , x))$ 
	\end{enumtag}
	Denote for convenience $\xi(t,e)$ as $t \centerdot e$; then \ref{ax_3} and \ref{ax_4} are the axioms
	\[\label{action_ofE}\forall t : e_0 \centerdot t = t \qquad\qquad \forall t,e,e' : (e \cdot e')\centerdot t = e\centerdot(e'\centerdot t)\]
	of a left action of $(E,\cdot,e_0)$ on $P$.

	In this notation, axioms \ref{ax_1} and \ref{ax_2} read
	\[
		\sigma (a , e \centerdot t) = \sigma (a \sact e , t)\qquad\qquad
		\delta (a , e \centerdot t) = (a \dact e) \centerdot (\delta (a \sact e , t))
	\]
	the first of which asserts the fact that $\sigma : A\times P \to X$ is `balanced' with respect to the structure of left $E$-sets on $(P,\centerdot),(A,\sact)$ (as a consequence, $\sigma$ descends to a map from the tensor product of $E$-sets $\hat\sigma : A\otimes_E P \to X$). As for the second, this is precisely the condition stated in \autoref{crucialemma} ensuring that the two actions of $A^*$ (induced by $\delta$) and $E$ (in \eqref{action_ofE}) on $P$ match to an action of $E\bowtie A^*$.
\end{savvyProof}
\begin{savvyProof}[Proof of \autoref{has_cumpa}]
	We build the companion of $f : A\to B$ as a cell
	\[\begin{tikzcd}[cramped, ampersand replacement=\&]
			{A\times 1} \& A \& B \& {1\times B}
			\arrow["\sim", from=1-1, to=1-2]
			\arrow["f", from=1-2, to=1-3]
			\arrow["\sim", from=1-3, to=1-4]
		\end{tikzcd}\]
	for which the `unit' and `counit' cells, boiling down to the diagrams
	\[\begin{tikzcd}[cramped, ampersand replacement=\&]
			{A\times 1} \& A \& B \& {1\times B} \& {A\times 1} \&\&\& {1\times A} \\
			{B\times 1} \&\&\& {1\times B} \& {A\times 1} \& A \& B \& {1\times B,}
			\arrow["\sim", from=1-1, to=1-2]
			\arrow["{f\times 1}"', from=1-1, to=2-1]
			\arrow["f", from=1-2, to=1-3]
			\arrow["\sim", from=1-3, to=1-4]
			\arrow[Rightarrow, no head, from=1-4, to=2-4]
			\arrow["\shortmid"{marking}, Rightarrow, no head, from=1-5, to=1-8]
			\arrow[Rightarrow, no head, from=1-5, to=2-5]
			\arrow["{1\times f}", from=1-8, to=2-8]
			\arrow["\shortmid"{marking}, Rightarrow, no head, from=2-1, to=2-4]
			\arrow["\sim"', from=2-5, to=2-6]
			\arrow["f"', from=2-6, to=2-7]
			\arrow["\sim"', from=2-7, to=2-8]
		\end{tikzcd}\]
	are readily constructed; the zig and zag identities are also totally elementary.
\end{savvyProof}
\begin{proof}[Proof of \autoref{crucialemma}]
	From \autoref{sbriggosbraggo}, $H\bowtie G$ comes equipped with a universal choice of a cospan $H \xto{i_H}H\bowtie G\xot{i_G} G$; precomposing a left representation $H\bowtie G \to\Set$ with $i_H, i_G$ induces a pair of left representations $H \to\Set, G\to\Set$ subject to condition \eqref{gndgoadbgoadb}, as it can be immediately seen from the relation $(1,g)\bullet (h,1) = (g\dact h,g\sact h)$ present in $H\bowtie G$.

	Conversely, given two representations $\alpha,\beta$ satisfying \eqref{gndgoadbgoadb}, we just have to show that posing $(h,g)\bullet x := \alpha(h, \beta(g,x))$ defines an action of $H\bowtie G$. The fact that $(1,1)\star x = x$ is immediate. As for associativity, one easily finds that
	\begin{align*}
		\big((h,g)\bullet (h',g')\big)\star x & = (h . g\dact h , g \sact h . g') \star x                      \\
		                                      & =\alpha(h . g\dact h , \beta(g \sact h . g',x) )               \\
		                                      & =\alpha(h, \alpha(g\dact h , \beta(g \sact h \beta(g',x) ) ) ) \\
		                                      & =\alpha(h, \beta(g , \alpha(h', \beta(g',x))))                 \\
		                                      & =\alpha(h, \beta(g , (h',g')\star x))                          \\
		                                      & =(h,g) \star ((h',g')\star x)
	\end{align*}
	which concludes the proof.
\end{proof}

\end{document}